\documentclass{article}
\usepackage[T2A]{fontenc}
\usepackage[cp1251]{inputenc}
\usepackage[english,russian]{babel}
\usepackage[tbtags]{amsmath}
\usepackage{amsfonts,amssymb,mathrsfs,amscd}

\usepackage{cite}

\usepackage[hyper]{msb-a}

\makeatletter
\gdef\No{{\select@language{russian}\textnumero}}
\makeatother

\JournalName{}
\numberwithin{equation}{section}
\theoremstyle{plain}
\newtheorem{theorem}{Теорема}

\newtheorem{propos}{Предложение}
\theoremstyle{definition}
\newtheorem{definition}{Определение}
\newtheorem{proof}{Доказательство}
\newtheorem{remark}{Замечание}

\begin{document}

\title{Operator-differential expressions: regularization and completeness of the root functions}
\author[S.\,A.~Buterin]{Sergey~Buterin}
\address{Saratov State University, Saratov, 410012, Russia;\\
Moscow Center for Fundamental and Applied Mathematics, Moscow, 119991, Russia;\\
Lomonosov Moscow State University, Moscow, 119991, Russia} \email{buterinsa@sgu.ru}

\udk{517.984}

\maketitle

\begin{fulltext}

\begin{abstract}
{\bf Abstract.} We consider an operator-differential expression of the form
$$
\ell y=\frac{d^m}{dx^m}\Big(By^{(n)}+Cy\Big), \quad 0<x<1,
$$
where $B$ is a linear bounded invertible operator, while $C$ is some finite-dimensi\-onal linear operator relatively bounded to the operator of $n$-fold differentiation. To such a form, we can reduce, in
particular, various singular differential expressions with the coefficients in negative Sobolev spaces, which creates an alternative to their regularization. In the case when $B$ is an integral Volterra
operator of the second kind with a continuous kernel vanishing at the diagonal, we establish completeness of the root functions of an operator generated by the expression $\ell y$ and irregular semi-separated
boundary conditions.


\end{abstract}

\begin{center}
{\large\bf Операторно-дифференциальные выражения: регуляризация и полнота корневых функций}
\end{center}

\smallskip
\begin{center}
С.\,А.~Бутерин
\end{center}

\address{Саратовский национальный исследовательский университет имени Н.Г. Чернышевского;\\
Московский центр фундаментальной и прикладной математики;\\
Московский государственный университет имени М.В. Ломоносова} \email{buterinsa@sgu.ru}

\bigskip
\begin{abstract}
Рассматривается операторно-дифференциальное выражение вида
$$
\ell y=\frac{d^m}{dx^m}\Big(By^{(n)}+Cy\Big), \quad 0<x<1,
$$
где $B$ -- линейный ограниченный обратимый оператор, а $C$ -- некоторый конечномерный линейный оператор, ограниченно подчиненный оператору $n$-кратного дифференцирования. Установлено, что к такому виду
приводятся, в частности, различные сингулярные дифференциальные выражения с коэффициентами из негативных соболевских пространств, что создает альтернативу их регуляризации. Для случая, когда $B$ является
интегральным оператором Вольтерра второго рода с непрерывным ядром, тождественно равным нулю на диагонали, доказана полнота системы собственных и присоединенных функций оператора, порожденного выражением~$\ell
y$ и нерегулярными полураспадающимися краевыми условиями.
\end{abstract}


\markright{Operator-differential expressions}


Работа посвящена 75-летнему юбилею члена-корреспондента РАН профессора Московского университета А.А.~Шкаликова и 90-летнему юбилею
заслужен\-ного деятеля науки РФ профессора Саратовского университета А.П.~Хромова.

\section{Введение}\label{s1}

В настоящей работе вводятся специальные классы операторно-диффе\-ренци\-альных выражений, в виде которых, помимо прочего, можно представить сингулярные дифференциальные выражения с коэффициентами из негативных
соболевских пространств (пространств обобщенных функций). В частности, это создает альтернативу регуляризации Мирзоева--Шкаликова таких дифференциальных выражений, которая, как известно, заключается в
построении так называемой согласованной матрицы из класса Шина--Зеттла. Кроме того, в работе получено описание всех согласованных матриц, что позволяет выражать краевые условия в терминах любых возможных
наборов квазипроизводных, а также обосновывается невозможность расширения существующих классов коэффициентов-распределений. Также реализован общий подход к определению весовых негативных соболевских
пространств, основанный на расширении классического пространства основных (тестовых) функций.

На порождаемые введенными выражениями операторы распространяется теорема Шкаликова о полноте собственных и присоединенных функций (с.п.ф.) дифференциальных операторов с суммируемыми коэффициентами и
нерегулярными распадающимися краевыми условиями. В качестве следствия получена полнота с.п.ф. соответствующих сингулярных дифференциальных операторов, для которых вопросы полноты еще не изучались.

Полнота устанавливается путем сведения к теореме Хромова о полноте с.п.ф. конечномерных возмущений интегральных вольтерровых операторов, доказательство которой, в свою очередь, также использует прием Шкаликова,
основанный на тригонометрической выпуклости индикатора целой функции.

\subsection{}
Недавно в связи с исследованием задач оптимального управления для некоторых классов нелокальных операторов \cite{MZ-24, Kras-24} автор
пришел к рассмотрению операторно-дифференциального выражения вида
\begin{equation}\label{eq1.1}
\ell_0 y:=\frac{d^n}{dx^n}By^{(n)}(x), \quad 0<x<1, \quad n\in{\mathbb N},
\end{equation}
где $B$ -- некоторый линейный ограниченный оператор, взаимно однозначно отображающий пространство $L_2(0,1)$ на себя.

Например, в \cite{Kras-24} этот оператор имел вид $B=(I+K^*)(I+K),$ где $I$ -- тождественный оператор, $K$ -- интегральный оператор
Вольтерра
$$
Kf(x)=\int_0^x K(x,t)f(t)\,dt, \quad \int_0^1dx\int_0^x |K(x,t)|^2\,dt <\infty,
$$
а $K^*$ является оператором сопряженным к $K.$

Выражение (\ref{eq1.1}) естественно рассматривать на функциях из класса $W_2^n[0,1],$ причем внешнюю производную можно понимать как обычным
образом, так и в смысле теории распределений.

{\it Первое} подразумевает рассмотрение только таких $y\in W_2^n[0,1],$ для которых $z:=By^{(n)} \in W_1^n[0,1].$ Нетрудно показать, что
они (как и те $y,$ для которых $z$ является многочленом) образуют всюду плотное подмножество $L_2(0,1).$ При этом функции $y^{\langle
n+j\rangle}:=z^{(j)},$ $j=\overline{0,n},$ выполняют роль квазипроизводных, а выражение (\ref{eq1.1}) приобретает смысл только при
$y^{\langle n+j\rangle}\in W_1^1[0,1],$ $j=\overline{0,n-1}.$

{\it Второе} означает
$$
(\ell_0 y,\varphi)=(-1)^n\int_0^1By^{(n)}(x)\varphi^{(n)}(x)\,dx
$$
для всякой бесконечно дифференцируемой функции $\varphi$ с носителем в $(0,1),$ т.е. $\varphi\in C_0^\infty(0,1),$ где $(\ell_0 y,\varphi)$ является значением функционала $\ell_0 y$ на функции~$\varphi.$ В этом
случае выражение (\ref{eq1.1}) определено уже {\it для всех} $y\in W_2^n[0,1].$

Заметим, что оба способа понимания выражения (\ref{eq1.1}) совпадут, как только $\ell_0y\in L(0,1),$ т.е. если $y$ является решением
уравнения $\ell_0 y=f$ с суммируемой правой частью и, в частности, решением однородного уравнения $\ell_0 y=\lambda y.$

Выражения вида (\ref{eq1.1}) с различными операторами $B$ могут появляться при рассмотрении билинейных или полуторалинейных форм, возникающих в вариационных методах. При этом положительная определенность такого
оператора~$B$ влечет коэрцитивность соответствующей формы. Например, в \cite{Skub-97} в связи с задачей об успокоении систем управления с последействием рассматривался случай $n=1,$ когда $B$ является
разностным оператором.

Таким образом, вид (\ref{eq1.1}) удобен тем, что объединяет в себе различные конкретные выражения. При этом можно отказаться от самосопряженности оператора~$B.$ Более универсальным представляется выражение
\begin{equation}\label{eq1.2}
\ell y:=\frac{d^m}{dx^m}\Big(By^{(n)}+Cy\Big),
\end{equation}
где $B$ -- линейный ограниченный обратимый оператор в некотором функциональном банаховом пространстве, а $C$ -- конечномерный оператор вида
\begin{equation}\label{eq1.3}
Cy=\sum_{\nu=0}^{n-1} y^{(\nu)}(0)u_\nu(x), \quad u_\nu\in {\rm Im}B, \;\; \nu=\overline{0,n-1},
\end{equation}
причем  линейная независимость функций $u_\nu(x)$ не подразумевается.

В виде (\ref{eq1.2}), (\ref{eq1.3}) можно представить также различные операторы дробного дифференцирования. Например, пусть $C=0,$ а $B$
является оператором дробного интегрирования Римана--Лиувилля порядка $\alpha\in(0,1).$ Тогда $\ell y$ будет дробной производной
Римана--Лиувилля функции $y$ порядка $m-\alpha,$ если $n=0,$ или ее дробной производной Герасимова--Капуто порядка $n-\alpha,$ если $m=0.$

\subsection{}
Выражение (\ref{eq1.2}) интересно также тем, что, как будет показано ниже, в таком в виде можно представить сингулярные дифференциальные выражения с коэффициентами из негативных соболевских пространств
\cite{Kurzweil, Guggenheimer, Pfaff, White, Weidmann, Filippov, Derr, NZadeSchk,b3, Sav-01, SavShk03, HrynMyk-03, HrynMyk-03-1, HrynMyk-04-1, HrynMyk-04, HrynMyk-04-2, Vlad-04, SavShk-05, HrynMyk-06-0,
HrynMyk-06, SavShk-08, FreIgnYur-08, Ign-10, GorMikh, SavShk-10, Sad-10, Hryn-11, HrynProns-12, EGNT-13, EGNT-13-1, SavShk-14, Mirz-14, EGNT-14, EGNST-15, MirShk-16, Vlad17-arXiv, KonMirShk-18, MirShk19-arXiv,
KonMir-19, Gul-19, UguBai-20, Bond-21, Bond-21-1, Bond-22, ZhAoBo-22, Bond-22-1, Bond-22-2, KonMirShk-23, Kuz-23, Bond-23, Bond-23-1, Chit-23, Bond-24, Chit-25}. При этом, чтобы <<уместить>> в виде
(\ref{eq1.2}) выражение наибольшей степени сингулярности, величина $|m-n|$ должна быть, соответственно, наименьшей.

Систематическое изучение спектральных свойств и других вопросов для сингулярных дифференциальных операторов с коэффициентами из пространств обобщенных функций началось с работ Шкаликова совместно с его
учениками: Нейманом-Заде (см.~\cite{NZadeSchk}) и Савчуком (см.~\cite{b3}), посвященных, в частности, оператору Штурма--Лиувилля с потенциалом $q\in W_2^{-1}.$

Регуляризация дифференциальных выражений произвольных порядков в так называемой нормальной форме с коэффициентами из классов максимально допустимой сингулярности получена в \cite{MirShk-16} для выражений
четного порядка и в \cite{MirShk19-arXiv} -- для нечетного. Под регуляризацией понимается приведение к эквивалентному выражению в пространстве вектор-функций, но с регулярными (суммируемыми либо локально
суммируемыми) коэффициентами.

Получаемое здесь представление сингулярных дифференциальных выражений в виде (\ref{eq1.2}) может служить альтернативой их регуляризации.

Если старший коэффициент такого выражения равен единице, то соответствующий оператор $B$ в (\ref{eq1.2}) можно представить в виде
\begin{equation}\label{eq1.4}
B=I+K, \quad Kf(x)=\int_0^x K(x,t)f(t)\,dt.
\end{equation}
При этом выражениям четного порядка \cite{MirShk-16} отвечает ядро $K(x,t),$ удовлетворяющее неравенству $|K(x,t)|\le a(x)+a(t)$ с некоторой функцией $a\in L_2(0,1),$ и функции $u_\nu\in L_2(0,1),$
$\nu=\overline{0,n-1},$ в (\ref{eq1.3}). Выражениям нечетного порядка~\cite{MirShk19-arXiv} при $m=n+1$ будет соответствовать оценка $|K(x,t)|\le a(t)\in L(0,1),$ а также $u_\nu\in AC[0,1],$
$\nu=\overline{0,n-1}.$ Для выражений четного порядка с произвольным старшим коэффициентом потребуется рассмотрение оператора $B$ в весовом пространстве. В общем случае вид $B$ получен в теоремах~\ref{th6}
и~\ref{th7}.

\subsection{}
Другой целью настоящей работы является доказательство полноты системы с.п.ф. оператора $L,$ порожденного выражением $\ell y$ вида (\ref{eq1.2})--(\ref{eq1.4}) с непрерывным ядром $K(x,t),$ равным нулю на
диагонали, и некоторыми нелокальными краевыми условиями, включающими, в частности, произвольные нерегулярные распадающиеся краевые условия для дифференциального оператора с суммируемыми коэффициентами (см.
условия (\ref{eq2.2}) в разделе~\ref{s2.1}).

Внешним аналогом краевых условий (\ref{eq2.2}) в случае выражения (\ref{eq1.2}) естественно считать линейно независимые условия вида
\begin{equation}\label{eq1.5}
\sum_{k=0}^{N-1}\alpha_{j,k}y^{\langle k\rangle}(0)=0, \quad j=\overline{1,l},
\end{equation}
\begin{equation}\label{eq1.6}
\sum_{k=0}^{N-1}\beta_{j,k}y^{\langle k\rangle}(1)=0, \quad j=\overline{l+1,N},
\end{equation}
где $N=m+n,$ а $\alpha_{j,k}$ и $\beta_{j,k}$ -- комплексные числа, причем
\begin{equation}\label{eq1.7}
y^{\langle k\rangle}:=y^{(k)}, \;\; k=\overline{0,n-1}, \quad y^{\langle
k\rangle}:=\frac{d^{k-n}}{dx^{k-n}}\Big(By^{(n)}+Cy\Big), \;\; k=\overline{n,N-1}.
\end{equation}

Функции $y^{\langle k\rangle},$ определенные в (\ref{eq1.7}) при $ k=\overline{n,N-1},$ будем называть ({\it нелокальными}) {\it квазипроизводными} -- в отличие от обычных квазипроизводных, участвующих в
регуляризации сингулярных дифференциальных выражений с коэффициентами-распределениями и являющихся локальными.

Не ограничивая общности, условия (\ref{eq1.5}) и (\ref{eq1.6}) можно считать нормированными:
\begin{equation}\label{eq1.8}
\sum_{k=0}^{\sigma_j}\beta_{j,k}y^{\langle k\rangle}(0)=0, \quad j=\overline{1,l}, \quad \sum_{k=0}^{\sigma_j}\beta_{j,k}y^{\langle k\rangle}(1)=0, \quad j=\overline{l+1,N},
\end{equation}
где $\beta_{j,\sigma_j}=1$ при $j=\overline{1,N},$ а также $0\le\sigma_{1}<\ldots<\sigma_l\le N-1$ и
\begin{equation}\label{eq1.9}
0\le\sigma_{l+1}<\ldots<\sigma_N\le N-1.
\end{equation}

Таким образом, условия в (\ref{eq1.8}) при $j=\overline{l+1,N},$  равно как и (\ref{eq1.6}), являются нелокальными, коль скоро ядро $K(x,t)$ в (\ref{eq1.4}) ненулевое. Однако вместо них мы будем рассматривать
более общие нелокальные условия вида
\begin{equation}\label{eq1.10}
U_j(y):=y^{\langle \sigma_j\rangle}(1)+\int_0^1 y^{\langle \sigma_j\rangle}(t)\,d\Phi_j(t)+\sum_{k=0}^{N-1}\alpha_{j,k}y^{\langle k\rangle}(0)=0, \;\; j=\overline{l+1,N},
\end{equation}
где функции $\Phi_j$ принадлежат классу $BV[0,1]$ и непрерывны в точках~$0$ и~$1.$

Аналогично теореме~5 в \cite[Гл.~VI, \S~6]{KolmFom} можно показать, что если снабдить слагаемое $y^{\langle \sigma_j\rangle}(1)$ в (\ref{eq1.10}) числовым коэффициентом и позволить ему обращаться в нуль, то в
случае $\alpha_{j,k}=0,$ $k=\overline{\sigma_j+1,N-1},$ линейная форма $U_j(y)$ будет представлять собой общий вид линейного ограниченного функционала в банаховом пространстве
$C^{\langle\sigma_j\rangle}:=\{y:y^{\langle\nu\rangle}\in C[0,1], \nu=\overline{0,\sigma_j}\}$ с нормой
$$
\|y\|_{C^{\langle\sigma_j\rangle}}:=\max_{\nu=\overline{0,\sigma_j}}\|y^{\langle\nu\rangle}\|_{C[0,1]}.
$$
Таким образом, (\ref{eq1.10}) является наиболее общей формой записи краевых условий (\ref{eq1.6}), причем представимость последних в виде (\ref{eq1.10}) следует из доказательства теоремы~\ref{th10}. В нем также
будет установлено, что именно к виду (\ref{eq1.5}), (\ref{eq1.10}) с $\Phi_j\in AC[0,1]$ приводятся произвольные распадающиеся, а также некоторые полураспадающиеся краевые условия для рассмотренных
дифференциальных операторов в терминах любых локальных квазипроизводных, включая условия (\ref{eq2.2}), если записать соответствующее выражение в виде (\ref{eq1.2}).

В точке~$0$ мы будем вынуждены ограничиться рассмотрением локальных условий (\ref{eq1.5}), а также потребовать $l>N-l,$ что продиктовано используемым подходом к получению нашей теоремы о полноте с.п.ф.
оператора $L.$

\begin{theorem}\label{th1}
Пусть операторы $B$ и $C$ в (\ref{eq1.2}) имеют вид (\ref{eq1.4}) и (\ref{eq1.3}), соответственно, причем $u_\nu\in L(0,1),$
$\nu=\overline{0,n-1},$ а $K(x,t)$ непрерывно, и
\begin{equation}\label{eq1.11}
K(x,x)\equiv0.
\end{equation}
Тогда система всех с.п.ф. краевой задачи для уравнения
\begin{equation}\label{eq1.12}
\ell y(x)=\lambda y(x), \quad 0<x<1,
\end{equation}
с условиями (\ref{eq1.5}) и (\ref{eq1.10}), отвечающими требованиям (\ref{eq1.9}) и  $l>N-l>0,$ полна в пространстве $L_2(0,1).$
\end{theorem}

Поскольку равенство (\ref{eq1.12}) влечет абсолютную непрерывность квазипроизводных $y^{\langle k\rangle}$ при $k=\overline{n,N-1},$ все краевые условия заданы корректно.

Если в виде (\ref{eq1.2})--(\ref{eq1.4}) представить то или иное дифференциальное выражение, то условие (\ref{eq1.11}) будет соответствовать обращению в нуль коэффициента при $(N-1)$-й производной в этом
выражении, а непрерывность ядра $K(x,t)$ будет иметь место, например, в случае непрерывности первых регулярных первообразных остальных его коэффициентов (см. замечания~\ref{rem3} и~\ref{rem5}).

Теорема~\ref{th1} обобщает теорему Шкаликова о полноте с.п.ф. дифференциальных операторов с нерегулярными распадающимися краевыми условиями~\cite{Shk-76}. Однако ее доказательство основано на сведении к
результатам Хромова о полноте для конечномерных возмущений интегральных вольтерровых операторов~\cite{Khrom-04}, в доказательстве которой, в свою очередь, также используется прием из \cite{Shk-76}, апеллирующий
к тригонометрической выпуклости индикатора.

\medskip
 Истории последнего вопроса посвящен следующий параграф.
В~\S~\ref{s3} приводятся сведения о сингулярных дифференциальных выражениях с коэффициен\-тами-распределениями. В~\S~\ref{s4} дано описание всех матриц из класса Шина--Зет\-тла, порождающих одно и то же
квазидифференциальное выражение и получены формулы для пересчета различных наборов квазипроизводных друг через друга. Последнее позволит не зависеть от выбора квазипроизводных при записи краевых условий для
дифференциальных операторов с коэффициен\-тами-рас\-пределениями.
 В~\S~\ref{s5} дается определение весовых негативных соболевских пространств в полной шкале индекса суммируемости путем
надлежащего расширения класса основных функций, а также обосновывается невозможность расширения классов коэффициентов-рас\-пределений.
 В~\S~\ref{s6} показано, что рассматри\-ваемые сингулярные дифференциальные выражения приводятся к опера\-торно-диффе\-рен\-циаль\-но\-му виду (\ref{eq1.2}),
а в~\S~\ref{s7} для них строится некоторый набор локальных квазипроизводных $y^{[k]}$ путем выделения локальных составляющих в $y^{\langle k\rangle}.$ При этом получены явные формулы пересчета $y^{\langle
k\rangle}$ через $y^{[k]}$ и обратно. Кроме того, дано описание всех согласованных матриц для одного и того же сингулярного дифференциального выражения и получена полная параметризация различных наборов
квазипроизводных.
 Доказательство теоремы~\ref{th1} приводится в~\S~\ref{s8}.
 В заключительном \S~\ref{s9} в качестве следствия устанавливается теорема о полноте с.п.ф. сингулярных дифференциальных операторов с нерегулярными полураспадающимися краевыми условиями.

\section{Теоремы Шкаликова и Хромова о полноте}\label{s2}

\subsection{Дифференциальные операторы}\label{s2.1}
Рассмотрим краевую задачу для дифференциального уравнения
\begin{equation}\label{eq2.1}
y^{(N)}+p_{N-2}(x)y^{(N-2)}+\ldots+p_0(x)y=\lambda y, \quad 0<x<1,
\end{equation}
с суммируемыми коэффициентами $p_j(x),$ $j=\overline{0,N-2},$ и линейно независимыми распадающимися краевыми условиями
\begin{equation}\label{eq2.2}
\sum_{k=0}^{N-1}\alpha_{j,k}y^{(k)}(0)=0, \quad j=\overline{1,l}, \qquad \displaystyle\sum_{k=0}^{N-1}\alpha_{j,k}y^{(k)}(1)=0, \quad
j=\overline{l+1,N}.
\end{equation}
В случае $2l=N$ условия (\ref{eq2.2}) являются условиями типа Штурма, которые регулярны по Биркгофу и даже усиленно регулярны \cite{Neum}. В частности, это влечет убывание функции Грина по всем направлениям в
комплексной плоскости спектрального параметра~$\lambda$ и позволяет установить безусловную базисность в $L_2(0,1)$ системы с.п.ф. краевой задачи (\ref{eq2.1}), (\ref{eq2.2}) (см. также \cite{Shk-79, Shk-82}).

Однако если $2l\ne N,$ то условия (\ref{eq2.2}) нерегулярны, и функция Грина имеет экспоненциальный рост, который может наблюдаться по любому направлению плоскости параметра $\lambda.$ Чтобы исключить задачу
Коши, без ущерба для общности можно потребовать $l>N-l>0.$ В этом случае вопросы разложения по с.п.ф. и их полноты становятся значительно более трудными.

В 1951 году Келдышем \cite{Keld-51} была анонсирована полнота системы с.п.ф. краевой задачи (\ref{eq2.1}), (\ref{eq2.2}). Ее доказательство в случае аналитических коэффициентов, указанное Хромовым
\cite{Khrom-73-d}, использовало результат Сахновича \cite{Sakh-58} о существовании <<треугольного>> оператора преобразования. Однако, как следовало из работы Мацаева \cite{Mats-60}, требование аналитичности
было существенным.

В общем случае суммируемых коэффициентов доказательство полноты было дано Шкаликовым \cite{Shk-76} в 1976-м году, т.е. спустя четверть века
после выхода работы Келдыша \cite{Keld-51}. Для получения нужных оценок Шкаликов предложил оригинальный прием, основанный на
тригонометрической выпуклости индикатора целой функции. В результате им была доказана следующая теорема.

\begin{theorem}[(Шкаликов)]\label{th2}
Пусть $l>N-l>0.$ Тогда с.п.ф. краевой задачи (\ref{eq2.1}), (\ref{eq2.2}) образуют полную систему в $L_2(0,1).$
\end{theorem}

Различные вопросы, связанные с разложением по с.п.ф., включая их полноту, восходят к задаче обоснования метода Фурье разделения переменных.
Исследования в этом направлении получили большое распространение после работ Биркгофа \cite{Birk-1, Birk-2}, который выделил класс так
называемых регулярных двухточечных краевых условий, для которых ряды по с.п.ф. ведут себя в смысле сходимости, как и обычные
тригонометрические ряды Фурье \cite{Tamarkin-17, Stone, Tamarkin-27}.

Например, как сказано выше, в случае $2l=N$ с.п.ф. краевой задачи (\ref{eq2.1}), (\ref{eq2.2}) образуют безусловный базис в $L_2(0,1),$ и, стало быть, соответствующие разложения сходятся в среднем квадратичном.
Если же разлагаемая функция принадлежит $W_1^N[0,1]$ и удовлетворяет краевым условиям (\ref{eq2.2}), то согласно теореме~4 в \cite[\S\,5]{Neum} будет иметь место уже равномерная сходимость.

Однако, как показал Джексон \cite{Jackson}, всякая функция, разлагаемая в равномерно сходящийся ряд по с.п.ф. простейшего
нерегулярного оператора с распадающимися краевыми условиями
$$
y''', \quad y(0)=y'(0)=y(1)=0,
$$
является бесконечно дифференцируемой. Позже Гопкинс \cite{Hopkins} уточнил, что такая функция является аналитической, причем
каждая третья ее производная, начиная с нулевой, удовлетворяет краевым условиям в нуле.

Как следует из \cite{Freiling-86}, к тому же результату приведет сходимость в среднем квадратичном, а значит с.п.ф. не могут быть
базисом Шаудера в $L_2(0,1).$

Метод Гопкинса впоследствии применялся к спектральному анализу оператора (\ref{eq2.1}), (\ref{eq2.2}), но лишь в случае, когда коэффициенты
$p_j(x)$ аналитичны (см., например, \cite{Ward, Seifert, Eberhard}, а также более подробную библиографию в \cite{Neum, Khrom-04,
Freiling}).

В общем случае решение задачи о разложении по с.п.ф. оператора (\ref{eq2.1}), (\ref{eq2.2}) было получено Хромовым \cite{Khromov-63, Khromov-66, Khromov-76}. Для случая непрерывных коэффициентов в уравнении
(\ref{eq2.1}) им было установлено, что в равномерно сходящиеся ряды по с.п.ф. разлагаются опе\-ра\-торно-ана\-лити\-ческие функции Фаге \cite{Fage} и только они.

В \cite{Freiling-86} показано, что операторную аналитичность разлагаемой функции также влечет сходимость соответствующего ряда в
среднем квадратичном.

\subsection{Конечномерные возмущения вольтерровых операторов}
Также Хромовым был выделен важный класс операторов, вообще говоря, в абстрактном банаховом пространстве, для которого полученные им результаты о разложении по с.п.ф. допускают естественные обобщения
\cite{Khrom-04, Khrom-73-d, Khrom-73, Khrom-74}.

В том числе речь идет об интегральных операторах вида
\begin{equation}\label{eq2.3}
Af=Mf+\sum_{k=1}^d g_k(x)\int_0^1 f(t)v_k(t)\,dt, \quad Mf:=\int_0^x M(x,t)f(t)\,dt,
\end{equation}
к которому приводятся обратные к дифференциальным и интегро-дифферен\-циальным операторам Вольтерра с произвольными краевыми условиями.

Для операторов вида (\ref{eq2.3}) была также получена следующая теорема о полноте (теорема~4.6 в \cite{Khrom-04}), из которой вытекает и
утверждение теоремы~\ref{th2}.

\begin{theorem}[(Хромов)]\label{th3}
Пусть функции $M(x,t),$ $g_k(x)$ и $v_k(x)$ при $k=\overline{1,d}$ непрерывны и для некоторого $N>2d$ удовлетворяют следующим условиям:

1) существует такое целое число $p\in[0,N],$ что функции
$$
\frac{\partial^{i+j}}{\partial x^i\partial t^j}M(x,t), \quad i=\overline{0,p}, \quad j=\overline{0,N-p},
$$
непрерывны при $0\le x-t\le\delta$ для некоторого $\delta>0;$

2) выполняется асимптотика
$$
M(x,t)=\frac{(x-t)^{N-1}}{(N-1)!} + o((x-t)^N), \quad x-t\to +0;
$$

3) для некоторых целых $\chi_k,\varkappa_k\in[0,N-1]$ имеют место асимптотики
$$
\left.\begin{array}{c}
g_k(x)=x^{\chi_k}(1+o(1)), \quad x\to0,\\[3mm]
v_k(x)=(1-x)^{\varkappa_k}(a_k+o(1)), \quad x\to 1,
\end{array}\right\} \quad a_k\ne0,
\quad k=\overline{1,d},
$$
причем $\chi_j\ne\chi_k$ и $\varkappa_j\ne\varkappa_k$ для $j\ne k.$

Тогда система всех с.п.ф. оператора (\ref{eq2.3}) полна в $L_2(0,1).$
\end{theorem}

При этом Хромовым была установлена глубокая связь полноты с.п.ф. оператора (\ref{eq2.3}) с вопросом о порождающих элементах оператора $M,$ в то время как вопросы полноты относятся к разряду более трудных. А
именно, если нуль не является собственным значением оператора $A,$ и выполняются условия~2) и~3) теоремы~\ref{th3}, то для полноты с.п.ф. необходимо и достаточно, чтобы система функций $g_k$ была порождающей
для оператора~$M$ (теорема~4.3 в \cite{Khrom-04}).

Последнее, помимо указанных условий на~$g_k$ и условия~2) теоремы~\ref{th3}, обеспечивается ее условием~1), что, в свою очередь, вытекает из теоремы~4.4 в \cite{Khrom-04}, при доказательстве которой также был
использован оригинальный прием Шкаликова \cite{Shk-76}, связанный с тригонометрической выпуклостью индикатора.

Однако, как показал Мацнев~\cite{Mats-1}, если ослабить условие~1) теоремы~\ref{th3}, заменив в нем значение $N$ на $N_1<N,$ то можно построить пример нециклического оператора $M,$ т.е. не имеющего ни одной
порождающей функции. Согласно сказанному выше это означает, что система с.п.ф. соответствующего оператора (\ref{eq2.3}) при $d=1$ никогда не будет полной в $L_2(0,1),$ что, в свою очередь, иллюстрирует
чрезвычайную тонкость приведенных результатов о полноте.

В то же время теорема о разложении (теорема~3.4 в \cite{Khrom-04}) данного условия~1) не содержит, и ее условия являются менее
ограничительными, чем условия теоремы о полноте. Поэтому вопрос полноты требовал отдельного исследования.

\section{Сингулярные дифференциальные выражения с коэффициентами-распределениями}\label{s3}

\subsection{Второй порядок}\label{s3.1}
Работа \cite{b3} положила начало систематическому исследованию спектральных свойств операторов, порождаемых выражением
\begin{equation}\label{eq3.1}
ly:=-y''+ qy
\end{equation}
с потенциалом $q\in W_2^{-1}[0,1]$ (в том числе см. \cite{Sav-01, SavShk03, HrynMyk-03, HrynMyk-03-1, HrynMyk-04-1, HrynMyk-04,
HrynMyk-04-2, SavShk-05, HrynMyk-06-0, HrynMyk-06, SavShk-08, FreIgnYur-08, Ign-10, GorMikh, SavShk-10, Sad-10, Hryn-11, HrynProns-12,
EGNT-13, EGNT-13-1, SavShk-14, Mirz-14, EGNT-14, EGNST-15, Gul-19, Bond-21, Kuz-23, Chit-23, Chit-25}).

В частности, в цикле работ \cite{SavShk-08, SavShk-10, Hryn-11, SavShk-14}, охватывающем для потенциала~$q$ полную шкалу гильбертовых пространств $W_2^{\theta}[0,1]$ при $\theta\ge-1,$  получены результаты по
обратным спектральным задачам, идея которых является новой даже для классического случая. Кроме того, основные результаты \cite{SavShk-10, SavShk-14} стали возможными в том числе за счет рассмотрения $q\in
W_2^\theta[0,1]$ при~$\theta<0.$

Принадлежность $q$ классу $W_2^{-1}[0,1]$ означает
\begin{equation}\label{eq3.2}
(q,\varphi)=-\int_0^1\sigma(x)\varphi'(x)\,dx
\end{equation}
для всех $\varphi\in C_0^\infty(0,1)$ с некоторой фиксированной функцией $\sigma\in L_2(0,1).$

Если функция $y$ является бесконечно дифференцируемой, то произведение $qy$ также принадлежит $W_2^{-1}[0,1]$ и в рамках классической теории обобщенных функций (см., например, \cite{KolmFom}) определяется
формулой
\begin{equation}\label{eq3.3}
(qy,\varphi):=(q,\varphi y).
\end{equation}
Однако в этом случае и $-y''$ будет бесконечно дифференцируемой функцией, а значит, не сможет выполнять роль главной части выражения (\ref{eq3.1}). Кроме того, соответствующее однородное уравнение в классе
таких $y$ может иметь, вообще говоря, только тривиальное решение.

Поэтому выражение (\ref{eq3.1}) необходимо было задать для всех $y\in W_2^1[0,1].$ При этом для определения $qy$ в \cite{NZadeSchk} была
построена теория мультипликаторов. А именно, {\it мультипликатором} в данном случае называется такая обобщенная функция $q,$ умножение на
которую представляет собой ограниченный оператор, действующий из $W_2^1[0,1]$ в $W_2^{-1}[0,1].$ Замечательный факт состоит в том, что~$q$
является мультипликатором тогда и только тогда, когда $q\in W_2^{-1}[0,1].$

Однако для определения произведения $qy$ можно поступить иначе, а именно, расширить пространство основных функций $C_0^\infty(0,1)$ вплоть
до ${\stackrel{\circ}{W}}\vphantom{W}_{\!2}^1[0,1].$

В самом деле, на $W_2^{-1}[0,1]$ можно смотреть как на пространство линейных ограниченных функционалов, действующих на
${\stackrel{\circ}{W}}\vphantom{W}_{\!2}^1[0,1].$ При этом функционал $q\in W_2^{-1}[0,1]$ также имеет вид (\ref{eq3.2}) с
некоторой функцией $\sigma\in L_2(0,1),$ заданной с точностью до константы, но теперь уже для $\varphi\in
{\stackrel{\circ}{W}}\vphantom{W}_{\!2}^1[0,1].$

В то же время согласно теореме Рисса об общем виде линейного ограниченного функционала в гильбертовом пространстве правая часть
(\ref{eq3.2}) представляет собой общий вид такого функционала на ${\stackrel{\circ}{W}}\vphantom{W}_{\!2}^1[0,1],$ поскольку в
(\ref{eq3.2}) можно ограничиться только теми $\sigma\in L_2(0,1),$ которые имеют нулевое среднее.

Таким образом, формула (\ref{eq3.2}) порождает изометрический изоморфизм между пространством $W_2^{-1}[0,1]$ с нормой $\|q\|_{W_2^{-1}[0,1]}:=\|\sigma\|_{L_2(0,1)}$ и пространством
${\stackrel{\circ}{W}}\vphantom{W}_{\!2}^1[0,1]$ с нормой $\|g\|_{{\stackrel{\circ}{W}}\vphantom{W}_{\!2}^1[0,1]}:=\|\sigma\|_{L_2(0,1)},$ причем в последнем случае $\sigma=g'.$

При таком определении $W_2^{-1}[0,1]$ произведение функционала $q\in W_2^{-1}[0,1]$ на функцию $y\in W_2^1[0,1]$ определяется формулой (\ref{eq3.3}) непосредственно и порождает ограниченный оператор,
переводящий $y\in W_2^1[0,1]$ в $qy\in W_2^{-1}[0,1],$ т.е. мультипликатор в смысле \cite{NZadeSchk} (см.~предложение~\ref{pr4}). В~\S~\ref{s5} этот подход распространяется на весовые негативные соболевские
пространства.

\subsection{Регуляризация}\label{s3.2}
Чтобы применять методы теории дифференциальных и интегральных уравнений, выражение (\ref{eq3.1}) необходимо регуляризовать. Другими словами, его нужно преобразовать к эквивалентному выражению в пространстве
вектор-функций, но с суммируемыми коэффициентами.

Для этой цели в \cite{b3} использовалась квазипроизводная
\begin{equation}\label{eq3.4}
y^{[1]}=y'-\sigma y,
\end{equation}
приводящая к соответствующему квазидифференциальному выражению
\begin{equation}\label{eq3.5}
-(y^{[1]})'-\sigma y^{[1]}-\sigma^2y.
\end{equation}
Похожие конструкции встречались и ранее (см. \cite{Guggenheimer, Pfaff}, а также обзор в \cite{KonMir-19}).

Выражение (\ref{eq3.5}) порождает тот же функционал из $W_2^{-1}[0,1],$ что и выражение (\ref{eq3.1}), а в случае $q\in L(0,1)$ и $y\in W_1^2[0,1]$ совпадает с (\ref{eq3.1}) почти всюду.

Согласно (\ref{eq3.4}) и (\ref{eq3.5}) выражение (\ref{eq3.1}) преобразуется к эквивалентному выражению
\begin{equation}\label{eq3.6}
\left[\begin{array}{c}
0\\[0mm]
-ly
\end{array}\right]=Y'-{\bf Q}Y, \quad {\bf Q}:=\left[\begin{array}{cc}
\sigma & 1\\[0mm]
-\sigma^2 & -\sigma
\end{array}\right],
\end{equation}
в пространстве вектор-функций $Y=[y,y^{[1]}]^T,$ но с суммируемой матрицей ${\bf Q}.$ При этом свободный член для (\ref{eq3.6}) должен иметь вид $[0,f]^T.$

В соответствии с \cite{MirShk-16} матрица ${\bf Q}$ является согласованной с выражением $-ly$ и принадлежит в терминологии \cite{EverMark-99} классу Шина--Зеттла (см. также раздел~\ref{s3.3}). Своя
согласованная матрица из этого класса возникает при любом способе регуляризации и хранит информацию о форме использованной квазипроизводной (либо квазипроизводных, если речь идет о выражениях высших порядков).

Как показала Бондаренко \cite{Bond-23}, все согласованные с выражением $-ly$ матрицы имеют вид
\begin{equation}\label{eq3.7}
\left[\begin{array}{cc}
\sigma_1 & 1\\[0mm]
q_1-\sigma_1^2 & -\sigma_1
\end{array}\right],
\quad \sigma_1\in L_2(0,1), \quad  q_1\in L(0,1), \quad q=\sigma_1' +q_1.
\end{equation}

Например, если для какой-либо цели требуется, чтобы $q_1=\sigma_1^2,$ то последнее соотношение в (\ref{eq3.7}) становится уравнением Риккати относительно $\sigma_1.$ В частности, его разрешимость в $L_2(0,1)$
при вещественных $q\in W_2^{-1}[0,1]$ равносильна тому, что $q$ является так называемым потенциалом Миуры (см. \cite{MZ-24}).

Различные формы квазипроизводных служат также для выражения краевых условий, заданных в точках, когда обычные производные не являются
непрерывными функциями и по этой причине непригодны для этой цели.

Зафиксируем некоторую функцию $q_1\in L(0,1)$ и обозначим через $y^{\{1\}}$ соответствующую квазипроизводную, т.е. $y^{\{1\}}=y'-\sigma_1y$
с какой-либо функцией $\sigma_1\in L_2(0,1),$ удовлетворяющей последнему соотношению в (\ref{eq3.7}).

Вычитая $y^{\{1\}}$ из (\ref{eq3.4}), получим $y^{[1]}-y^{\{1\}}=(\sigma_1-\sigma)y.$ С другой стороны, согласно определению
$\sigma$ и $\sigma_1,$ будем иметь $\sigma'-\sigma_1'=q_1\in L(0,1).$ В результате приходим к следующему утверждению.

\begin{propos}\label{pr1}
Все (локальные) квазипроизводные для выражения~$-ly$ образуют однопараметрическое семейство $y^{[1]}+fy,$ подчиненное функциональному параметру $f\in AC[0,1],$ где $y^{[1]}$ -- определенная, например, в
(\ref{eq3.4}).
\end{propos}

Таким образом, локальные квазипроизводные для выражения (\ref{eq3.1}) отличаются друг от друга слагаемым вида $fy,$ где $f\in AC[0,1].$ Данное обстоятельство позволяет не зависеть от выбора конкретной
квазипроизводной при задании краевых условий к выражению (\ref{eq3.1}). Другими словами, если краевые условия содержат некоторую квазипроизводную, то их можно эквивалентным образом выразить через какую-либо
другую квазипроизводную. При этом любые две квазипроизводные для выражения (\ref{eq3.1}) могут быть абсолютно непрерывными только одновременно. В~\S~\ref{s4} аналогичные свойства устанавливаются для общего
квазидифференциального выражения произвольного порядка, а в~\S~\ref{s7} -- распространяются на произвольные дифференциальные выражения с коэффициентами-распределениями в нормальной форме.

\subsection{Квазидифференциальное выражение, согласованная матрица}\label{s3.3}
В \cite{MirShk-16} для дифференциального выражения некоторого общего вида с коэффици\-ентами-распределениями было введено понятие
согласованной матрицы, которое без изменения распространяется на любое дифференциальное выражение. В дальнейшем нам понадобится
это понятие применительно к выражениям, заданным на интервале $(0,1).$ Чтобы привести соответствующее определение, напомним общую
конструкцию квазидифференциального выражения, предложенную Шином \cite{Shin} и позже использованную Зеттлом~\cite{Zettl-75} и
другими.

Рассмотрим матри\-цу-функ\-цию ${\bf F}=[f_{k,j}]_{k,j=1}^N,$ у которой все элементы $f_{k,j}$ принадлежат $L(0,1),$ причем
$f_{k,k+1}\ne0$ п.в. на $(0,1)$ при $k=\overline{1,N-1}$ и $f_{k,j}=0$ при $j>k+1.$ При этих условиях будем говорить, что ${\bf
F}$ принадлежит классу ${\cal S}_N[0,1],$ являющемуся, в свою очередь, подклассом класса ${\cal S}_N(0,1)$ матриц Шина--Зеттла с
локально суммируемыми элементами (см. \cite{MirShk-16}).

Матрица ${\bf F}$ порождает квазидифференциальное выражение
\begin{equation}\label{eq3.8}
{\cal F}y: =(y^{[N-1]})'-\sum_{j=1}^Nf_{N,j}y^{[j-1]},
\end{equation}
где квазипроизводные $y^{[k]}$ рекуррентно определяются формулами
\begin{equation}\label{eq3.9}
y^{[0]}:=y, \quad y^{[k]}:=\frac1{f_{k,k+1}}\Big((y^{[k-1]})'-\sum_{j=1}^kf_{k,j}y^{[j-1]}\Big), \quad k=\overline{1,N-1}.
\end{equation}
Множество ${\cal D}({\cal F}):=\{y:y^{[k]}\in AC[0,1],\, k=\overline{0,N-1}\}$ называется областью определения выражения (оператора) ${\cal F}y.$

Очевидно, принадлежность $y\in {\cal D}({\cal F})$ влечет ${\cal F}y\in L(0,1).$ В силу теоремы~1 в \cite[\S\,16]{Neum} справедливо и обратное. Точнее, для всякой функции $f\in L(0,1)$ и любого числа
$x_0\in[0,1]$ существует и единственно решение уравнения ${\cal F}y=f$ с наперед заданными значениями $y^{[k]}(x_0),$ $k=\overline{0,N-1}.$

Рассмотрим какое-либо дифференциальное выражение $\ell y,$ заданное на интервале $(0,1).$ Следующее определение соответствует определению~2 в \cite{MirShk-16}.

\begin{definition}\label{def0}
Матрицу ${\bf F}\in {\cal S}_N[0,1]$ назовем {\it согласованной} с выражением $\ell y,$ если оно может быть корректно определено как распределение для всех $y\in{\cal D}({\cal F})$ и совпадет на них с ${\cal
F}y,$ т.е. $\ell y={\cal F}y$ в смысле распределений.
\end{definition}

Последнее равенство будет выполняться и почти всюду, если~$\ell y$ считать суммируемой функцией, что становится возможным в силу суммируемости ${\cal F}y.$

Например, тот факт, что матрица в (\ref{eq3.7}) является согласованной с выраже\-нием $-ly$ означает, что порождаемое ею квазидифференциальное выражение
$$
(y'-\sigma_1y)'-(q_1-\sigma_1^2)y +\sigma_1(y'-\sigma_1y)
$$
совпадает с $-ly$ на тех $y\in AC[0,1],$ для которых $y'-\sigma_1 y\in AC[0,1].$ При этом все такие $y,$ очевидно, принадлежат
классу $W_2^1[0,1],$ в точности на котором выражение $ly$ может быть определено как распределение.

Вообще, {\it регуляризация} какого-либо сингулярного дифференциального выражения $\ell y$ есть не что иное, как построение согласованной матрицы~${\bf F}$ из соответствующего  класса Шина--Зеттла.

\subsection{Дифференциальные выражения произвольного порядка}\label{s3.4}
В \cite{MirShk-16} и \cite{MirShk19-arXiv} получена регуляризация сингулярных дифференциальных выражений четного и нечетного порядков, соответственно, в так называемой нормальной форме (см. также \cite{Vlad-04,
Vlad17-arXiv, Bond-22, Bond-23}). В \cite{MirShk-16} рассматривалось выражение
\begin{equation}\label{eq3.10}
\begin{array}{c}
\displaystyle \ell_{2n} y:=\sum_{k=0}^n (p_k(x)y^{(n-k)})^{(n-k)} \qquad\qquad\qquad\qquad\qquad\qquad\qquad\quad\\[5mm]
\displaystyle \qquad\qquad\quad+i\sum_{k=0}^{n-1} \Big((q_k(x)y^{(n-k-1)})^{(n-k)} +(q_k(x)y^{(n-k)})^{(n-k-1)}\Big),
\end{array}
\end{equation}
коэффициенты которого удовлетворяют условиям
\begin{equation}\label{eq3.11}
\sqrt{|p_0|}, \frac1{\sqrt{|p_0|}}, \frac{P_k}{\sqrt{|p_0|}}, \frac{Q_{k-1}}{\sqrt{|p_0|}}\in L_2(0,1), \quad k=\overline{1,n},
\end{equation}
где $P_k$ и $Q_k$ являются $k$-ми первообразными $p_k$ и $q_k,$ соответственно, в смысле распределений. В свою очередь, в \cite{MirShk19-arXiv} было рассмотрено выражение вида
\begin{equation}\label{eq3.12}
\begin{array}{c}
\displaystyle \ell_{2n+1}y:=i\sum_{k=0}^{n} \Big((q_{k}(x)y^{(n-k+1)})^{(n-k)} +(q_{k}(x)y^{(n-k)})^{(n-k+1)}\Big) \qquad\qquad\\[5mm]
\displaystyle \qquad\qquad\qquad \qquad\qquad\qquad\qquad\qquad\qquad+ \sum_{k=0}^n (p_k(x)y^{(n-k)})^{(n-k)}
\end{array}
\end{equation}
со следующими условиями на коэффициенты:
\begin{equation}\label{eq3.13}
q_0,\frac1{q_0}\in AC[0,1], \quad p_0, P_k, Q_k\in L(0,1), \quad k=\overline{1,n},
\end{equation}
где, как и выше, $P_k^{(k)}=p_k,$ но теперь $Q_k^{(k-1)}=q_k$ при $k=\overline{1,n}.$

В \cite{MirShk-16, MirShk19-arXiv} получены, в частности, согласованные матрицы ${\bf Q},$ позволяющие записать соответствующие
выражения $\ell_Ny$ в эквивалентном векторном виде
\begin{equation}\label{eq3.14}
{\cal L}Y=Y'-{\bf Q}Y,
\end{equation}
где $Y=[y^{[0]},\ldots,y^{[N-1]}]^T,$ ${\cal L}Y=[0,\ldots,0,\ell_Ny]^T,$ а $y^{[j]}$ являются соответствующими квазипроизводными, причем $y^{[j]}=y^{(j)}$ для $j=\overline{0,n-1}.$ Кроме того, свободный член
для (\ref{eq3.14}) должен иметь нули на тех же местах, что и~${\cal L}Y.$

Опираясь на результаты работ \cite{MirShk-16, MirShk19-arXiv}, Бондаренко \cite{Bond-21-1, Bond-22, Bond-22-1, Bond-22-2, Bond-23, Bond-23-1, Bond-24} построила общую теорию обратных спектральных задач для
дифференциальных операторов высших порядков с коэффициентами-распределениями.

Также Бондаренко \cite{Bond-23} получено специальное семейство согласованных матриц для выражений (\ref{eq3.10}) и (\ref{eq3.12}) со старшим коэффициентом единицей. В разделе~\ref{s7.4} будет указано, как
построить полное семейство в общем случае, а в разделе~\ref{s4.3} --  для произвольного квазидифференциального выражения.

Как и в разделе~\ref{s3.1}, выражения (\ref{eq3.10}) и (\ref{eq3.12}) можно определить в соответствующих пространствах распределений (замечание~\ref{rem1}). Однако если этого не делать, то квазипроизводные
могут послужить для придания им смысла. С такой целью в \cite{Krein, Glaz, Neum} использовалась более очевидная, но вместе с тем и более <<слабая>> форма квазипроизводных для частного случая (\ref{eq3.10}) на
произвольном интервале $(a,b),$ когда $\frac1{p_0},p_1,\ldots,p_n\in L_{\rm loc}(a,b),$ а все $q_k=0.$

В работах \cite{MirShk-16, MirShk19-arXiv} выражения (\ref{eq3.10}) и (\ref{eq3.12}) тоже рассматривались на произвольном интервале, вообще говоря, бесконечном, а пространства коэффициентов обладали индексом
<<loc>>. В связи с этим согласованные матрицы ${\bf Q}$ были локально суммируемы. Однако рассмотрение интервала~$(0,1)$ и пространств без индекса <<loc>> в отношении регуляризации общности не уменьшает.

Отметим также, что Курцвейлем~\cite{Kurzweil}
 рассматривались дифференциальные уравнения
произвольных порядков в полиномиальной форме (\ref{eq2.1}) (в том числе когда присутствует и слагаемое с $(N-1)$-й производной) с коэф\-фици\-ентами-рас\-пределениями (см. также \cite{Filippov}). Изучению
представления соответствующих выражений в виде (\ref{eq1.2}) планируется уделить внимание в отдельной работе.

\section{Полуэквивалентные и эквивалентные матрицы в ${\cal S}_N[0,1]$}\label{s4}

Все рассуждения и результаты настоящего параграфа ес\-тественным образом распространяются на классы локально суммируемых матриц
Шина--Зеттла ${\cal S}_N(a,b)$ при $-\infty\le a<b\le\infty$ (см. определение в \cite{MirShk-16}).

\subsection{Определение}\label{s4.1}
Тогда как ${\bf F}\in {\cal S}_N[0,1]$ определяет соответствующее квазидифференциальное выражение (\ref{eq3.8}) однозначно, обратное неверно, т.е. найдется другая матрица в ${\cal S}_N[0,1],$ порождающая то же
выражение ${\cal F}y.$

Рассматривая помимо~${\bf F}$ матрицу $\widetilde{\bf F}=[\widetilde f_{k,j}]_{k,j=1}^N\in{\cal S}_N[0,1],$ будем обозначать
через $\widetilde{\cal F}y$ квазидифференциальное выражение, порождаемое последней, т.е.
$$
\widetilde{\cal F}y: =(y^{\{N-1\}})'-\sum_{j=1}^N\widetilde f_{N,j}y^{\{j-1\}},
$$
где
\begin{equation}\label{eq4.1}
y^{\{0\}}:=y, \quad y^{\{k\}}:=\frac1{\widetilde f_{k,k+1}}\Big((y^{\{k-1\}})'-\sum_{j=1}^k\widetilde f_{k,j}y^{\{j-1\}}\Big), \quad k=\overline{1,N-1}.
\end{equation}

\begin{definition}\label{def1} Матрицы ${\bf F},\widetilde{\bf F}\in {\cal S}_N[0,1]$ будем называть {\it полуэквивалентными}, если имеет место
\begin{equation}\label{eq4.2}
{\cal D}({\cal F})={\cal D}(\widetilde{\cal F}),
\end{equation}
и {\it эквивалентными}, если они полуэквивалентны и справедливо тождество
\begin{equation}\label{eq4.3}
{\cal F}y=\widetilde{\cal F}y \quad \forall y\in{\cal D}({\cal F}).
\end{equation}
\end{definition}

\subsection{Полуэквивалентные матрицы}\label{s4.2}
Следующая теорема позволяет параметризовать наборы квазипроизводных, порождаемые всеми матрицами из класса ${\cal S}_N[0,1],$ являющимися полуэквивалентными между собой.

\begin{theorem}\label{th4}
1) Для любого набора функций
\begin{equation}\label{eq4.4}
g_{k,j},\frac1{g_{k,k}}\in AC[0,1], \quad k=\overline{1,N-1}, \quad j=\overline{0,k},
\end{equation}
найдется такая матрица $\widetilde{\bf F}\in{\cal S}_N[0,1],$ что будут справедливы соотношения
\begin{equation}\label{eq4.5}
y^{\{k\}}=\sum_{j=0}^kg_{k,j}y^{[j]}, \quad k=\overline{1,N-1}, \quad \forall y\in {\cal D}({\cal F}),
\end{equation}
а значит,
\begin{equation}\label{eq4.6}
{\cal D}({\cal F})\subset{\cal D}(\widetilde{\cal F}).
\end{equation}

2) Обратно, если матрица $\widetilde{\bf F}\in{\cal S}_N[0,1]$ такова, что выполняется (\ref{eq4.6}), то существует набор функций (\ref{eq4.4}), для которого имеет место (\ref{eq4.5}).
\end{theorem}

\begin{proof} 1) Согласно (\ref{eq4.1}) требуется найти такие функции
\begin{equation}\label{eq4.7}
\widetilde f_{k,j}\in L(0,1), \quad k=\overline{1,N-1}, \quad j=\overline{1,k+1},
\end{equation}
чтобы для левых частей (\ref{eq4.5}) выполнялось
\begin{equation}\label{eq4.8}
(y^{\{k-1\}})'=\sum_{j=1}^{k+1}\widetilde f_{k,j}y^{\{j-1\}}, \quad k=\overline{1,N-1},  \quad \forall y\in {\cal D}({\cal F}).
\end{equation}
Подставляя (\ref{eq4.5}) в (\ref{eq4.8}), для каждого $k=\overline{1,N-1}$ будем иметь
$$
\sum_{j=0}^{k-1}\Big(g_{k-1,j}'y^{[j]} +g_{k-1,j}(y^{[j]})'\Big)=\sum_{j=1}^{k+1}\widetilde f_{k,j}\sum_{\nu=0}^{j-1}g_{j-1,\nu}y^{[\nu]}, \quad g_{0,0}:=1.
$$
Выражая здесь $(y^{[j]})'$ при помощи (\ref{eq3.9}), получаем
$$
\sum_{j=0}^{k-1}g_{k-1,j}'y^{[j]} + \sum_{j=0}^{k-1}g_{k-1,j}\sum_{\nu=1}^{j+2} f_{j+1,\nu}y^{[\nu-1]}=\sum_{j=1}^{k+1}\widetilde f_{k,j}\sum_{\nu=0}^{j-1}g_{j-1,\nu}y^{[\nu]}.
$$
Меняя порядок суммирования, имеем
\begin{equation}\label{eq4.9}
\sum_{j=0}^k a_{k,j}y^{[j]}=0, \quad k=\overline{1,N-1},  \quad \forall y\in {\cal D}({\cal F}),
\end{equation}
где
\begin{equation}\label{eq4.10}
a_{k,j}=g_{k-1,j}' +\sum_{\nu=j+\delta_{0,j}}^kg_{k-1,\nu-1} f_{\nu,j+1} -\sum_{\nu=j}^kg_{\nu,j} \widetilde f_{k,\nu+1}, \quad j=\overline{0,k-1},
\end{equation}
\begin{equation}\label{eq4.11}
a_{k,k}=g_{k-1,k-1}f_{k,k+1}-g_{k,k}\widetilde f_{k,k+1}.
\end{equation}

Пусть ${\bf Y}:=[y_\nu^{[j-1]}]_{j,\nu=1}^N$ является фундаментальной матрицей системы дифференциальных уравнений
\begin{equation}\label{eq4.12}
(y^{[k-1]})'=\sum_{j=1}^{k+1-\delta_{k,N}}f_{k,j}y^{[j-1]}, \quad k=\overline{1,N}.
\end{equation}
Поскольку последняя равносильна уравнению ${\cal F}y=0,$ все $y_\nu$ принадлежат множеству ${\cal D}({\cal F}).$ Подставляя их в (\ref{eq4.9}) вместо $y$ и учитывая, что $\det{\bf Y}$ не обращается в нуль,
получаем $a_{k,j}=0$ при $k=\overline{1,N-1}$ и $j=\overline{0,k}.$

Таким образом, (\ref{eq4.11}) непосредственно дает формулы
\begin{equation}\label{eq4.13}
\widetilde f_{k,k+1}=\frac{g_{k-1,k-1}}{g_{k,k}}f_{k,k+1}, \quad k=\overline{1,N-1},
\end{equation}
а (\ref{eq4.10}) приводит к линейным алгебраическим системам
\begin{equation}\label{eq4.14}
G_k\widetilde f_k=[h_{k,1},\ldots,h_{k,k}]^T \in (L(0,1))^k, \quad k=\overline{1,N-1},
\end{equation}
относительно вектор-функций $\widetilde f_k:=[\widetilde f_{k,1},\ldots,\widetilde f_{k,k}]^T,$ где
\begin{equation}\label{eq4.15}
G_k:=[g_{\nu-1,j-1}]_{j,\nu=1}^k, \quad g_{\nu,j}:=0, \quad \nu<j,
\end{equation}
$$
h_{k,j}=g_{k-1,j-1}'-g_{k,j-1}\widetilde f_{k,k+1}+\sum_{\nu=j-1+\delta_{1,j}}^kg_{k-1,\nu-1} f_{\nu,j}, \quad j=\overline{1,k}.
$$

Согласно (\ref{eq4.4}) матрица $G_k^{-1}$ существует и ее элементы принадлежат $AC[0,1].$ Итак, мы построили набор функций (\ref{eq4.7}), для которого выполняется (\ref{eq4.8}). Остается заметить, что
содержащиеся в этом наборе функции (\ref{eq4.13}) почти всюду отличны от нуля, что в итоге дает первое утверждение теоремы.

2) Пусть матрицы ${\bf F}$ и $\widetilde{\bf F}$ из ${\cal S}_N[0,1]$ таковы, что выполняется (\ref{eq4.6}). Другими словами, $[y^{[0]},\ldots,y^{[N-1]}]\in (AC[0,1])^N$ влечет
$[y^{\{0\}},\ldots,y^{\{N-1\}}]\in (AC[0,1])^N.$ Покажем, что имеют место соотношения (\ref{eq4.5}) с некоторыми функциями (\ref{eq4.4}).

Заметим, что (\ref{eq4.5}) будет тривиальным для $k=0,$ если положить ${g_{0,0}:=1.}$ Предположим по индукции, что соотношения в~(\ref{eq4.5}) справедливы при всех $k=\overline{0,s-1}$ для некоторого
$s\in\{1,\ldots,N-1\},$ причем участвующие в них $g_{k,j}$ удовлетворяют соответствующим условиям в~(\ref{eq4.4}). Тогда, согласно пункту~1) доказательства, (\ref{eq4.5}) будет иметь место и для $k=s,$ т.е.
выполнится
\begin{equation}\label{eq4.16}
y^{\{s\}}=\sum_{j=0}^sg_{s,j}y^{[j]}  \quad \forall y\in {\cal D}(\widetilde{\cal F}),
\end{equation}
если $g_{s,j}$ при $j=\overline{0,s}$ выбрать таким образом, чтобы для $k=s$ были справедливы соотношения (\ref{eq4.13}) и (\ref{eq4.14}). Последнее означает, что нужно положить
\begin{equation}\label{eq4.17}
g_{s,s}=\frac{f_{s,s+1}}{\widetilde f_{s,s+1}}g_{s-1,s-1},
\end{equation}
\begin{equation}\label{eq4.18}
\begin{array}{c}
\displaystyle g_{s,j}=\frac1{\widetilde f_{s,s+1}}\Big(g_{s-1,j}'+\sum_{\nu=j-1+\delta_{0,j}}^{s-1}g_{s-1,\nu} f_{\nu+1,j+1}\qquad\qquad\qquad\qquad\qquad\\[3mm]
\displaystyle \qquad\qquad\qquad\qquad\qquad\qquad\qquad\qquad-\sum_{\nu=j+1}^s g_{\nu-1,j}\widetilde f_{s,\nu}\Big), \;\; j=\overline{0,s-1}.
\end{array}
\end{equation}
Таким образом, для завершения доказательства теоремы остается показать
\begin{equation}\label{eq4.19}
g_{s,j}\in AC[0,1], \quad j=\overline{0,s}.
\end{equation}

Рассмотрим фундаментальную матрицу ${\bf Y}=[y_\nu^{[j-1]}]_{j,\nu=1}^N$ системы (\ref{eq4.12}). Поскольку $y_\nu\in{\cal D}(\cal F),$ $\nu=\overline{1,N},$ то в силу (\ref{eq4.6}) имеем $y_\nu\in{\cal
D}(\widetilde{\cal F})$ и, в частности, $y_\nu^{\{s\}}\in AC[0,1]$ при $\nu=\overline{1,N}.$ Тогда из (\ref{eq4.16}) вытекают соотношения
\begin{equation}\label{eq4.20}
y_\nu^{\{s\}}=\sum_{j=0}^sg_{s,j}y_\nu^{[j]}, \quad \nu=\overline{1,N}.
\end{equation}
Так как существует обратная матрица ${\bf Y}^{-1}$ и все ее элементы, очевидно, принадлежат классу $AC[0,1],$ то соотношения (\ref{eq4.20}) дают
$$
[g_{s,0},\ldots,g_{s,s},0,\ldots,0]=[y_1^{\{s\}},\ldots,y_N^{\{s\}}]{\bf Y}^{-1},
$$
и мы приходим к (\ref{eq4.19}). $\hfill\Box$
\end{proof}

{\sc Следствие 1. }{\it Если матрицы ${\bf F},\widetilde{\bf F}\in {\cal S}_N[0,1]$ таковы, что имеет место включение (\ref{eq4.6}), то справедливо и обратное включение, а значит, выполняется (\ref{eq4.2}),
т.е. данные матрицы полуэквивалентны.}

\medskip
{\sc Доказательство.} Согласно второму утверждению теоремы~\ref{th4} существует такой набор функций (\ref{eq4.4}), что выполняется (\ref{eq4.5}). В то же время первое утверждение теоремы гарантирует для всякого
набора функций
\begin{equation}\label{eq4.21}
\widetilde g_{k,j},\frac1{\widetilde g_{k,k}}\in AC[0,1], \quad k=\overline{1,N-1}, \quad j=\overline{0,k},
\end{equation}
существование матрицы $\widetilde{\widetilde{\bf F}}=[{\,\,\widetilde{\!\!\widetilde f}_{k,j}}]_{k,j=1}^N\in{\cal S}_N[0,1],$ обеспечивающей
\begin{equation}\label{eq4.22}
y^{\{\{k\}\}}=\sum_{j=0}^k\widetilde g_{k,j}y^{\{j\}}, \quad k=\overline{1,N-1}, \quad \forall y\in {\cal D}(\widetilde{\cal F}).
\end{equation}
Здесь $y^{\{\{k\}\}}$ -- квазипроизводные, порожденные $\widetilde{\widetilde{\bf F}},$ т.е. выполняется
\begin{equation}\label{eq4.23}
(y^{\{\{k-1\}\}})'=\sum_{j=1}^{k+1}{\,\,\widetilde{\!\!\widetilde f}}_{k,j}y^{\{\{j-1\}\}}, \quad k=\overline{1,N-1},
\end{equation}
где $y^{\{\{0\}\}}=y.$ Таким образом, для квазидифференциального выражения ${\,\widetilde{\!\widetilde{\cal F}}}y,$ порожденного матрицей $\widetilde{\widetilde{\bf F}},$ имеет место включение
\begin{equation}\label{eq4.24}
{\cal D}(\widetilde{\cal F})\subset{\cal D}({\,\widetilde{\!\widetilde{\cal F}}}).
\end{equation}

Зададим функции (\ref{eq4.21}) формулой $\widetilde G_N=G_N^{-1},$ где матрица $G_N$ определена в (\ref{eq4.15}) при $k=N,$ а $\widetilde G_N$ строится аналогичным образом по $\widetilde g_{k,j}.$
 Тогда, поскольку соотношения (\ref{eq4.5}) и (\ref{eq4.22}) могут быть записаны в матричной форме:
$$
[y^{\{0\}},\ldots,y^{\{N-1\}}]=[y^{[0]},\ldots,y^{[N-1]}]G_N \quad \forall y\in {\cal D}({\cal F}),
$$
$$
[y^{\{\{0\}\}},\ldots,y^{\{\{N-1\}\}}]=[y^{\{0\}},\ldots,y^{\{N-1\}}]\widetilde G_N \quad \forall y\in {\cal D}(\widetilde{\cal F}),
$$
легко увидеть, что включение (\ref{eq4.6}) дает
$$
[y^{\{\{0\}\}},\ldots,y^{\{\{N-1\}\}}]=[y^{[0]},\ldots,y^{[N-1]}] \quad \forall y\in {\cal D}({\cal F}).
$$
Следовательно, вычитая (\ref{eq4.23}) из (\ref{eq4.12}) при $k=\overline{1,N-1},$ будем иметь
$$
\sum_{j=1}^{k+1}(f_{k,j}-{\,\,\widetilde{\!\!\widetilde f}}_{k,j})y^{[j-1]}=0, \quad k=\overline{1,N-1}, \quad \forall y\in {\cal D}({\cal F}).
$$
Подставляя сюда вместо $y^{[j-1]}$ элементы $y_\nu^{[j-1]}$ фундаментальной матрицы~${\bf Y}$ системы (\ref{eq4.12}), получаем ${\,\,\widetilde{\!\!\widetilde f}}_{k,j}=f_{k,j}$ при $k=\overline{1,N-1}$ и
$j=\overline{1,k+1}.$ Поскольку ${\cal D}({\cal F})$ вполне определяется первыми $N-1$ строками матрицы ${\bf F},$ приходим к равенству ${\cal D}({\,\widetilde{\!\widetilde{\cal F}}})={\cal D}({\cal F}),$
которое вместе с (\ref{eq4.24}) дает ${\cal D}(\widetilde{\cal F})\subset{\cal D}({\cal F}). $ $\hfill\Box$

\subsection{Эквивалентные матрицы}\label{s4.3}
Уточним связь между наборами квазипроизводных, порождаемыми теперь уже эквивалентными матрицами.

\begin{theorem}\label{th5}
Для всякой матрицы ${\bf F}\in{\cal S}_N[0,1]$ и любого набора функций (\ref{eq4.4}), в котором
\begin{equation}\label{eq4.25}
g_{N-1,N-1}\equiv1,
\end{equation}
существует единственная эквивалентная~${\bf F}$ матрица $\widetilde{\bf F}\in{\cal S}_N[0,1]$ такая, что выполняется (\ref{eq4.5}).

С другой стороны, если матрицы ${\bf F}, \widetilde{\bf F}\in{\cal S}_N[0,1]$ эквивалентны, то в соответствующем представлении (\ref{eq4.5}) имеет место (\ref{eq4.25}).
\end{theorem}

\begin{proof}
Существование матрицы $\widetilde{\bf F},$ для которой выполняется (\ref{eq4.5}), следует из первого утверждения теоремы~\ref{th4}. Кроме того, согласно его доказательству все ответственные за это первые ее
$N-1$ строк однозначно определяются заданием функций (\ref{eq4.4}). Далее, поскольку добавляемое условие (\ref{eq4.3}) является не чем иным, как равенством вида (\ref{eq4.5}) при $k=N,$
$y^{\{N\}}:=\widetilde{\cal F}y$ и $y^{[N]}:={\cal F}y,$ а также $g_{N,N}:=1$ и $g_{N,j}:=0,$ $j=\overline{0,N-1},$ остается повторить рассуждения доказательства первого утверждения теоремы~\ref{th4} для~$k=N.$

В самом деле, предполагая (\ref{eq4.3}), для $k=N$ аналогичным образом получаем соотношения (\ref{eq4.9})--(\ref{eq4.11}), где $f_{N,N+1}=\widetilde f_{N,N+1}=1.$ Пусть $[y_\nu^{[j-1]}]_{j,\nu=1}^{N+1}$
явля\-ется фундаментальной матрицей системы дифференциальных уравнений
$$
(y^{[k-1]})'=\sum_{j=1}^{k+1}f_{k,j}y^{[j-1]}, \quad k=\overline{1,N}, \quad (y^{[N]})'=0.
$$
Очевидно, все $y_\nu\in{\cal D}(\cal F).$ Подставляя их в (\ref{eq4.9}) при $k=N,$ получим $a_{N,j}=0,$ $j=\overline{0,N},$ что с учетом (\ref{eq4.11}) и $g_{N,N}\equiv1$ дает (\ref{eq4.25}), а с учетом
(\ref{eq4.10}) -- систему
\begin{equation}\label{eq4.26}
G_N[\widetilde f_{N,1},\ldots,\widetilde f_{N,N}]^T=[h_{N,1},\ldots,h_{N,N}]^T
\end{equation}
для нахождения последней строки матрицы $\widetilde{\bf F},$ где
$$
h_{N,j}=g_{N-1,j-1}'+\sum_{\nu=j-1+\delta_{1,j}}^Ng_{N-1,\nu-1} f_{\nu,j}, \quad j=\overline{1,N},
$$
а матрица $G_N$ определяется в (\ref{eq4.15}) при $k=N.$ $\hfill\Box$
\end{proof}

Итак, мы показали, что наборы квазипроизводных, порождаемые двумя эквивалентными матрицами, связаны соотношениями (\ref{eq4.5}) посредством некоторых абсолютно непрерывных функций вида (\ref{eq4.4}),
(\ref{eq4.25}). Последние, в свою очередь, могут быть найдены рекуррентно с помощью формул (\ref{eq4.17}) и (\ref{eq4.18}).

Теоремы~\ref{th4} и~\ref{th5} позволяют построить все матрицы $\widetilde{\bf F}$ эквивалентные произвольной матрице ${\bf F}.$ Для этого нужно рассмотреть всевозможные наборы абсолютно непрерывных функций вида
(\ref{eq4.4}), (\ref{eq4.25}). Тогда элементы соответствующих $\widetilde{\bf F}$ находятся по формулам (\ref{eq4.13}), где $g_{0,0}\equiv1,$ и из треугольных алгебраических систем (\ref{eq4.14}) и
(\ref{eq4.26}).

Теорема~\ref{th5} вместе со следствием~1 в предыдущем разделе и следствием~2 в разделе~\ref{s7.3} позволит получить обобщение предложения~\ref{pr1} на случай произвольного сингулярного выражения $\ell_Ny$ вида
(\ref{eq3.10}) или (\ref{eq3.12}) (теорема~\ref{th9}).

\section{Негативные соболевские пространства}\label{s5}

\subsection{Основные функции}
Пусть $1\le p\le\infty$ и $p^*=\frac{p}{p-1}.$ Рассмотрим почти всюду отличную от нуля функцию $b\in L_p(0,1)$ такую, что $\frac1b\in
L_{p*}(0,1).$

Обозначим через $L_{p,b}$ банахово пространство измеримых на $(0,1)$ функций~$f,$ удовлетворяющих условию $bf\in L_p(0,1),$ с нормой
$\|f\|_{p,b}:=\|bf\|_{L_p(0,1)}.$

Также через $W_{p,b}^k$ будем обозначать банахово пространство функций $g$ таких, что $g,g',\ldots,g^{(k-1)}\in AC[0,1]$ и $g^{(k)}\in L_{p,b},$ с нормой
$$
\|g\|_{W_{p,b}^k}:=\left\{\begin{array}{cc}
\displaystyle\sqrt{\sum_{\nu=0}^{k-1}|g^{(\nu)}(0)|^2+\|g^{(k)}\|_{2,b}^2}, & p=2,\\[5mm]
\displaystyle\max_{\nu=\overline{0,k-1}}|g^{(\nu)}(0)|+\|g^{(k)}\|_{p,b}, & p\ne2,
\end{array}\right.
$$
где $k+1\in{\mathbb N}.$ В частности, будем иметь $W_{p,b}^0=L_{p,b}.$

Отметим, что при $k\in{\mathbb N}$ определение $W_{p,b}^k$ корректно в силу включения $L_{p,b}\subset L(0,1),$ вытекающего при $1<p<\infty$
из неравенства Г\"ельдера:
$$
\int_0^1|f(x)|\,dx \le \Big(\int_0^1|f(x)b(x)|^p\,dx\Big)^\frac1p \Big(\int_0^1\frac{dx}{|b(x)|^{p^*}}\Big)^\frac1{p^*}<\infty
$$
и из соответствующих неравенств для норм при $p=1$ и $p=\infty.$

Следующее утверждение очевидно.

\begin{propos}\label{pr2}
Справедливы включения
$$
W_{p,b}^k\subset W_1^k[0,1], \quad k+1\in{\mathbb N},
$$
\begin{equation}\label{eq5.1}
W_{p,b}^k\subset W_{p,b}^{k-1}, \quad k\in{\mathbb N},
\end{equation}
\begin{equation}\label{eq5.2}
 fg\in W_{p,b}^k \quad \forall f,g\in W_{p,b}^k, \quad k\in{\mathbb N}.
\end{equation}
\end{propos}

При этом (\ref{eq5.1}) и (\ref{eq5.2}) вытекают из очевидного включения $L_\infty(0,1)\subset L_{p,b},$ которое позволит также рассмотреть замкнутое подпространство ${\stackrel{\circ}{W}}\vphantom{W}_{p,b}^k$
пространства $W_{p,b}^k,$ состоящее из функций, удовлетворяющих краевым условиям $g^{(\nu)}(0)=g^{(\nu)}(1)=0,$ $\nu=\overline{0,k-1}.$ Ниже оно будет играть роль пространства основных (тестовых) функций при
определении негативных пространств.

\subsection{Негативные пространства}\label{s5.2}
Обозначим через $W_{p^*\!,\frac1b}^{-k}$ множество функционалов $f$ на пространстве ${\stackrel{\circ}{W}}\vphantom{W}_{\!p,b}^k,$ имеющих вид
\begin{equation}\label{eq5.3}
(f,\varphi)=(-1)^k\int_0^1\sigma(x)\varphi^{(k)}(x)\,dx, \quad \varphi\in {\stackrel{\circ}{W}}\vphantom{W}_{\!p,b}^k,
\end{equation}
с некоторой функцией $\sigma\in L_{p^*\!,\frac1b}.$

Например, всякая функция $f\in L(0,1)$ задает такой функционал формулой
$$
(f,\varphi)=\int_0^1f(x)\varphi(x)\,dx, \quad \varphi\in {\stackrel{\circ}{W}}\vphantom{W}_{\!p,b}^k.
$$
Действительно, интегрируя $k$ раз по частям, приходим к (\ref{eq5.3}) с $\sigma=J^k f,$ где
\begin{equation}\label{eq5.4}
Jf:= \int_0^x f(t)\,dt.
\end{equation}
Поэтому функционал $f$ в (\ref{eq5.3}) естественно называть $k$-ой (обобщенной) производной функции $\sigma.$ При этом функционал $f$ определяет функцию $\sigma$ с точностью до многочлена степени меньшей $k$ и
является расширением с $C_0^\infty(0,1)$ на ${\stackrel{\circ}{W}}\vphantom{W}_{\!p,b}^k$ функционала, соответствующего ее $k$-ой обобщенной производной в смысле классической теории обобщенных функций (см.
\cite{KolmFom}).

Можно без потери общности считать, что выполняются условия
\begin{equation}\label{eq5.5}
\int_0^1(1-x)^\nu\sigma(x)\,dx=0, \quad \nu=\overline{0,k-1},
\end{equation}
и задать норму в $W_{p^*\!,\frac1b}^{-k}$ формулой $\|f\|_{W_{p^*\!,\frac1b}^{-k}}:=\|\sigma\|_{p^*\!,\frac1b}.$

Следующее утверждение вытекает из общего вида линейного ограниченного функционала на $L_p(0,1)$ при $p\in[1,\infty)$ (см., например,
\cite{LjustSob}).

\begin{propos}\label{pr3}
При $p\in[1,\infty)$ формула (\ref{eq5.3}) представляет собой общий вид линейного ограниченного функционала на ${\stackrel{\circ}{W}}\vphantom{W}_{\!p,b}^k.$
\end{propos}

Таким образом, $W_{p^*\!,\frac1b}^{-k}$ является реализацией пространства, сопряженного пространству
${\stackrel{\circ}{W}}\vphantom{W}_{\!p,b}^k,$ коль скоро $p\in[1,\infty).$ Для $p=\infty$ это не так, что связано с нерефлексивностью
пространства $L(0,1).$

Однако при всех $p\in[1,\infty]$ формула (\ref{eq5.3}) дает изометрический изоморфизм между $W_{p^*\!,\frac1b}^{-k}$ и ${\stackrel{\circ}{W}}\vphantom{W}_{\!p^*\!,\frac1b}^k,$ ставящий в соответствие
функционалу $f\in W_{p^*\!,\frac1b}^{-k}$ функцию $J^k\sigma\in {\stackrel{\circ}{W}}\vphantom{W}_{\!p^*\!,\frac1b}^k,$ где $\sigma$ определена в (\ref{eq5.3}) и удовлетворяет условиям (\ref{eq5.5}).

{\it Производной} $f\in W_{p^*\!,\frac1b}^{-k}$ назовем функционал $f',$ определяемый формулой
\begin{equation}\label{eq5.6}
(f',\varphi):=-(f,\varphi') \quad \forall\varphi\in {\stackrel{\circ}{W}}\vphantom{W}_{\!p,b}^{k+1}.
\end{equation}

Поскольку согласно (\ref{eq5.2}) множество $W_{p,b}^k$ является замкнутым относительно поточечного умножения функций, можно непосредственно определить произведение $fy$ функционала $f\in W_{p^*\!,\frac1b}^{-k}$
на функцию $y\in W_{p,b}^k$ формулой
\begin{equation}\label{eq5.7}
(fy,\varphi):=(f,\varphi y)  \quad \forall  \varphi\in{\stackrel{\circ}{W}}\vphantom{W}_{\!p,b}^k.
\end{equation}

Как и в случае выражения (\ref{eq3.1}), для определения сингулярных  выражений (\ref{eq3.10}) и (\ref{eq3.12}) в пространстве распределений ключевое значение имеет понятие мультипликатора. Здесь этому отвечает
следующее утверждение.

\begin{propos}\label{pr4}
Всякий функционал $f\in W_{p^*\!,\frac1b}^{-k}$ порождает линейный ограниченный оператор, переводящий $y\in W_{p,b}^k$ в $fy\in W_{p^*\!,\frac1b}^{-k}$ по формуле (\ref{eq5.7}).
\end{propos}

\begin{proof}
Пусть $f\in W_{p^*\!,\frac1b}^{-k}.$ Убедимся в справедливости оценки
\begin{equation}\label{eq5.8}
|(fy,\varphi)|\le C\|y\|_{W_{p,b}^k}\|\varphi\|_{W_{p,b}^k}  \quad \forall y\in W_{p,b}^k,  \quad \forall
\varphi\in{\stackrel{\circ}{W}}\vphantom{W}_{\!p,b}^k,
\end{equation}
где $C$ не зависит от $y$ и $\varphi.$ Согласно (\ref{eq5.3}) и (\ref{eq5.7}) будем иметь
$$
|(fy,\varphi)|\le\|\sigma\|_{p^*\!,\frac1b}\|(\varphi y)^{(k)}\|_{p,b},
$$
причем
$$
\|(\varphi y)^{(k)}\|_{p,b} \le \sum_{\nu=0}^k C_k^\nu \|\varphi^{(\nu)} y^{(k-\nu)}\|_{p,b} \le C_1\|\varphi\|\|y\|,
$$
где $C_1$ не зависит от $y$ и $\varphi,$ а норма
$$
\|y\|:=\max_{\nu=\overline{0,k-1}}\|y^{(\nu)}\|_{L_\infty(0,1)}+\|y^{(k)}\|_{p,b},
$$
эквивалентна норме $\|y\|_{W_{p,b}^k},$ что приводит к (\ref{eq5.8}). Таким образом, (\ref{eq5.7}) определяет на ${\stackrel{\circ}{W}}\vphantom{W}_{\!p,b}^k$ линейный ограниченный функционал $fy.$ При
$p\in[0,\infty)$ предложение~\ref{pr3} гарантирует $fy\in W_{p^*\!,\frac1b}^{-k},$ что можно показать и непосредственно, охватив также случай $p=\infty.$ Для этого, в соответствии с определением множества
$W_{p^*\!,\frac1b}^{-k},$ достаточно указать для каждого $y\in{W_{p,b}^k}$ такое $\sigma_y\in L_{p^*\!,\frac1b},$ что
$$
(fy,\varphi)=(-1)^k\int_0^1\sigma_y(x)\varphi^{(k)}(x)\,dx, \quad \varphi\in {\stackrel{\circ}{W}}\vphantom{W}_{\!p,b}^k.
$$
С учетом (\ref{eq5.3}) и (\ref{eq5.7}) нетрудно убедиться, что для этой цели подойдет
$$
\sigma_y=\sum_{\nu=0}^k C_k^\nu(-1)^\nu J^\nu(\sigma y^{(\nu)}),
$$
что завершает доказательство. $\hfill\Box$
\end{proof}

\begin{remark}\label{rem1}
Формула (\ref{eq5.7}) вместе с предложением~\ref{pr4} позволяет непосредственно определить сингулярные выражения (\ref{eq3.10}) и (\ref{eq3.12}) в пространстве распределений. А именно, выражение четного порядка
(\ref{eq3.10}) таким образом определено для всех $y\in W_{2,b}^n$ и представляет собой функционал из $W_{2,\frac1b}^{-n},$ где $b^2=p_0.$ В свою очередь, выражение нечетного порядка (\ref{eq3.12}) определено
для всех $y\in W_{\infty,1}^n$ и является функционалом из $W_{\infty,1}^{-n-1}.$
\end{remark}

В дальнейшем индекс $b$ в обозначениях введенных весовых пространств будем опускать, если $b=1.$ Например, $W_p^k$ будет обозначать
$W_{p,1}^k.$

\subsection{О невозможности расширения пространств сингулярных коэффициентов}
Рассмотрим выражение Штурма--Лиувилля (\ref{eq3.1}) с сингулярным потенциалом $q\in W_{p^*}^{-k}$ для некоторых $k\ge1$ и $p^*\ge1.$

В соответствии с разделом~\ref{s5.2} произведение $qy$ определено в точности для $y\in W_p^k.$ При этом $qy\in W_{p^*}^{-k}$ и $y''\in W_p^{k-2},$ а выражение (\ref{eq3.1}) понимается как функционал на
${\stackrel{\circ}{W}}\vphantom{W}_p^k$ при $k>1$ или на ${\stackrel{\circ}{W}}\vphantom{W}_{\max\{p,p^*\}}^1$ при $k=1.$

Поэтому если, например, $k>1,$ т.е. $k-2>-k,$ то второе слагаемое в (\ref{eq3.1}) становится <<более сингулярным>>, чем первое, которое по
этой причине не сможет выполнять роль главной части, как в случае $q\in W_2^{-1}.$

Если же $k=1,$ а $p^*<2,$ т.е. $p>2,$ то при одинаковой сингулярности обоих членов в (\ref{eq3.1}) второй снова будет принадлежать более широкому множеству, нежели первый, поскольку $W_p^{-1}$ является
собственным подмножеством $W_{p^*}^{-1}.$ В этом случае $-y''$ также не может оставаться главной частью.

Аналогичные попытки расширения класса какого-либо нестаршего коэффициента в выражениях (\ref{eq3.10}) и (\ref{eq3.12}) потребуют сужения класса функций $y,$ на которых получаемое выражение определено. В рамках
введенной шкалы пространств последнее приведет к тому, что член со старшей производной не сможет выполнять роль главной части. При этом соответствующее однородное уравнение будет иметь, вообще говоря, лишь
тривиальное решение.

\begin{remark}\label{rem2}
В \cite{ValNazSul-21, Bond-26} рассматривались некоторые выражения нечетного порядка вида (\ref{eq3.12}), коэффициенты которых не подчиняются условиям (\ref{eq3.13}). В частности, речь шла о выражении
\begin{equation}\label{eq5.9}
y'''+\sigma''y, \quad \sigma\in L_3(0,1).
\end{equation}
Например, в \cite{Bond-26} это выражение формально понималось, как квазидифференциальное выражение ${\cal F}y,$ порожденное матрицей
\begin{equation}\label{eq5.10}
\left[\begin{array}{ccc}
-\sigma & 1 & 0\\[0mm]
-\frac32\sigma^2 & 2\sigma & 1\\[0mm]
\sigma^3 & -\frac32\sigma^2 & -\sigma
\end{array}\right] \in {\cal S}_3[0,1].
\end{equation}
При этом матрица (\ref{eq5.10}), вообще говоря, не является согласованной с выражением (\ref{eq5.9}), поскольку оно определено в точности на функциях $y\in W_\frac32^2,$ не входящих в ${\cal D}({\cal F})$ (см.
замечание~2.2 в \cite{Bond-26}). В самом деле, $y^{[1]}(=y'+\sigma y)\in AC[0,1]$ может повлечь $y\notin W_1^2(\supset W_\frac32^2),$ коль скоро $\sigma\notin AC[0,1].$

Другими словами, связь (\ref{eq5.9}) и (\ref{eq5.10}) носит, вообще говоря, лишь формальный характер и состоит только в том, что (\ref{eq5.9}) совпадает с выражением ${\cal F}y$ при достаточно гладких $\sigma.$
В частности, при $\sigma\in AC[0,1]$ матрица (\ref{eq5.10}) становится согласованной с выражением (\ref{eq5.9}).
\end{remark}

\subsection{Свойства обобщенного дифференцирования}
Приведем три ут\-верждения, которые позволят в~\S~\ref{s6} использовать для преобразования выражений (\ref{eq3.10}) и (\ref{eq3.12}) <<обычные>> правила дифференцирования.

\begin{propos}\label{pr5}
Пусть $k+1\in{\mathbb N}.$ Для всех $f\in W_{p^*\!,\frac1b}^{-k}$ и $y\in W_{p,b}^{k+1}$ справедлива формула
\begin{equation}\label{eq5.11}
(fy)'=f'y+fy'.
\end{equation}
\end{propos}

\begin{proof}
Согласно (\ref{eq5.6}) и (\ref{eq5.7}) для всех $\varphi\in{\stackrel{\circ}{W}}\vphantom{W}_{\!p,b}^{k+1}$ имеем
$$
((fy)',\varphi)= -(f,\varphi'y) = (f,\varphi y'-(\varphi y)') =(f y',\varphi)+(f'y,\varphi),
$$
что и означает (\ref{eq5.11}). $\hfill\Box$
\end{proof}

Следующие два утверждения являются следствиями предложения~\ref{pr5}.

\begin{propos}\label{pr6}
Для всех $\sigma\in L_{p^*\!,\frac1b}$ и $y\in W_{p,b}^k,$ $k\in{\mathbb N},$ имеет место соотношение
\begin{equation}\label{eq5.12}
(\sigma y)^{(\nu)}=\sum_{j=0}^\nu C_\nu^j \sigma^{(j)}y^{(\nu-j)}, \quad \nu=\overline{1,k}.
\end{equation}
\end{propos}

\begin{proof}
При $k=1$ формула (\ref{eq5.12}) совпадает с (\ref{eq5.11}) для $k=0.$ Предположим по индукции, что (\ref{eq5.12}) имеет место при некотором $k\in{\mathbb N}$ и покажем ее справедливость для $k+1.$ Достаточно
рассмотреть $\nu=k+1.$ По предположению индукции, а также согласно предложению~\ref{pr5} будем иметь
$$
(\sigma y)^{(k+1)}=\sum_{j=0}^k C_k^j (\sigma^{(j)}y^{(k-j)})'=\sum_{j=0}^k C_k^j \sigma^{(j+1)}y^{(k-j)} +\sum_{j=0}^k C_k^j
\sigma^{(j)}y^{(k-j+1)},
$$
откуда после приведения подобных членов получаем (\ref{eq5.12}) для $\nu=k+1.$ $\hfill\Box$
\end{proof}

\begin{propos}\label{pr7}
Для всех $\sigma\in L_{p^*\!,\frac1b}$ и $y\in W_{p,b}^n,$ $n\in{\mathbb N},$ справедливо соотношение
\begin{equation}\label{eq5.13}
\sigma^{(k)}y^{(n-k)}=\sum_{\nu=0}^k (-1)^{\nu}C_k^\nu(\sigma y^{(n-k+\nu)})^{(k-\nu)}, \quad k=\overline{1,n}.
\end{equation}
\end{propos}

\begin{proof}
Можно прийти к формуле (\ref{eq5.13}), преобразуя ее правую часть при помощи~(\ref{eq5.12}), но проще воспользоваться индукцией по $k.$ В самом деле, для $k=0$ формула (\ref{eq5.13}) тривиальна. Предположим,
что она верна для некоторого $k\in\{0,\ldots,n-1\}.$ Тогда в силу предложения~\ref{pr5} будем иметь
$$
\sigma^{(k+1)}y^{(n-k-1)}=(\sigma^{(k)}y^{(n-k-1)})'-\sigma^{(k)}y^{(n-k)}
$$
$$
=\sum_{\nu=0}^k (-1)^{\nu}C_k^\nu(\sigma y^{(n-k+\nu-1)})^{(k-\nu+1)}-\sum_{\nu=0}^k (-1)^{\nu}C_k^\nu(\sigma y^{(n-k+\nu)})^{(k-\nu)}.
$$
Приводя подобные члены, приходим к (\ref{eq5.13}) для $k+1$ вместо $k.$ $\hfill\Box$
\end{proof}

\section{Приведение сингулярных дифференциальных выражений (\ref{eq3.10}) и (\ref{eq3.12}) к виду (\ref{eq1.2}), (\ref{eq1.3}) }\label{s6}

\subsection{Четный порядок}
Введем вес $b$ так, чтобы  $b^2=p_0.$ Тогда в соответствии с разделом~\ref{s5.2} условия (\ref{eq3.11}) примут вид
\begin{equation}\label{eq6.1}
p_0, \frac1{p_0}\in L(0,1), \quad p_k\in W_{2,\frac1b}^{-k}, \;\; k=\overline{1,n}, \quad q_\nu\in W_{2,\frac1b}^{-\nu}, \;\;\nu=\overline{0,n-1}.
\end{equation}

Согласно замечанию~\ref{rem1} выражение (\ref{eq3.10}) с коэффициентами, удовлетворяющими (\ref{eq6.1}), определено при $y\in W_{2,b}^n$ и является функционалом $\ell_{2n}y\in W_{2,\frac1b}^{-n}.$

\begin{theorem}\label{th6}
Для выражения (\ref{eq3.10}) с коэффициентами (\ref{eq6.1}) имеет место представление
\begin{equation}\label{eq6.2}
\ell_{2n}y=\frac{d^n}{dx^n}\Big(By^{(n)}+Cy\Big), \quad 0<x<1,
\end{equation}
где оператор $B$ взаимно однозначно отображает $L_{2,b}$ на $L_{2,\frac1b},$ причем
\begin{equation}\label{eq6.3}
Bf=b^2(x)f(x)+\int_0^x K(x,t)f(t)\,dt, \quad |K(x,t)|\le a(x)+a(t), \;\; a\in L_{2,\frac1b},
\end{equation}
\begin{equation}\label{eq6.4}
Cy=\sum_{\nu=0}^{n-1} y^{(\nu)}(0)u_\nu(x), \quad u_\nu(x)\in L_{2,\frac1b}, \quad \nu=\overline{0,n-1}.
\end{equation}
Кроме того, справедливо представление
$$
K(x,t)=\sum_{k=1}^n \Bigg(\frac{(x-t)^{k-1}}{(k-1)!}\Big(P_k(x)+iQ_{k-1}(x) +(-1)^k(P_k(t)-iQ_{k-1}(t))\Big)
$$
$$
+\sum_{\nu=1}^{k-1} (-1)^\nu\int_t^x\frac{(x-\tau)^{\nu-1}}{(\nu-1)!}\Big(C_k^\nu
P_k(\tau)+i(C_{k-1}^\nu-C_{k-1}^{\nu-1})Q_{k-1}(\tau)\Big) \frac{(\tau-t)^{k-1-\nu}}{(k-1-\nu)!}\,d\tau\Bigg),
$$
где $C_k^\nu$ -- биномиальные коэффициенты, $P_k^{(k)}=p_k$ и $Q_{k-1}^{(k-1)}=q_{k-1}$ при $k=\overline{1,n}.$
\end{theorem}

\begin{proof} Перепишем (\ref{eq3.10}) в виде
\begin{equation}\label{eq6.5}
\begin{array}{c}
\displaystyle \ell_{2n}y=\sum_{k=0}^n (P_k^{(k)}y^{(n-k)})^{(n-k)} \qquad\qquad\qquad\qquad\qquad\qquad\qquad\qquad\qquad\\[3mm]
\qquad\qquad\qquad \displaystyle +i\sum_{k=0}^{n-1} \Big((Q_k^{(k)}y^{(n-k-1)})^{(n-k)} +(Q_k^{(k)}y^{(n-k)})^{(n-k-1)}\Big),
\end{array}
\end{equation}
где $P_0=p_0.$ Поскольку $y\in W_{2,b}^n,$ в силу предложения~\ref{pr7} будем иметь
$$
P_k^{(k)}y^{(n-k)}=\sum_{\nu=0}^k (-1)^{\nu}C_k^\nu(P_k y^{(n-k+\nu)})^{(k-\nu)}, \quad k=\overline{0,n},
$$
$$
Q_k^{(k)}y^{(n-k-1)}=\sum_{\nu=0}^k (-1)^{\nu}C_k^\nu(Q_k y^{(n-k+\nu-1)})^{(k-\nu)}, \quad k=\overline{0,n-1},
$$
$$
Q_k^{(k)}y^{(n-k)}=\sum_{\nu=0}^k (-1)^{\nu}C_k^\nu(Q_k y^{(n-k+\nu)})^{(k-\nu)}, \quad k=\overline{0,n-1}.
$$
Подставив последние три представления в (\ref{eq6.5}), получаем
$$
\ell_{2n} y=\frac{d^n}{dx^n}\Bigg(\sum_{k=0}^n \sum_{\nu=0}^k (-1)^{\nu}C_k^\nu J^\nu(P_k y^{(n-k+\nu)})
\qquad\qquad\qquad\qquad\qquad
$$
\begin{equation}\label{eq6.6}
\qquad\qquad+i\sum_{k=0}^{n-1}\sum_{\nu=0}^k (-1)^{\nu}C_k^\nu  \Big(J^\nu(Q_k y^{(n-k+\nu-1)}) +J^{\nu+1}(Q_k y^{(n-k+\nu)})\Big)\Bigg),
\end{equation}
где оператор $J$ определен в (\ref{eq5.4}). В силу соотношения
\begin{equation}\label{eq6.7}
y^{(s)}=\sum_{j=s}^{n-1}y^{(j)}(0)\frac{x^{j-s}}{(j-s)!}+J^{n-s}y^{(n)}, \quad s=\overline{0,n},
\end{equation}
выражение (\ref{eq6.6}) имеет вид (\ref{eq6.2}), причем
\begin{equation}\label{eq6.8}
\begin{array}{c}
\displaystyle Bf=\sum_{k=0}^n \sum_{\nu=0}^k (-1)^{\nu}C_k^\nu J^\nu(P_k J^{k-\nu}f) \qquad\qquad\qquad\qquad\qquad\qquad\qquad\quad\\[3mm]
\displaystyle \quad\qquad\qquad +i\sum_{k=0}^{n-1} \sum_{\nu=0}^k (-1)^{\nu}C_k^\nu \Big(J^\nu(Q_k J^{k-\nu+1}f) +J^{\nu+1}(Q_k
J^{k-\nu}f)\Big),
\end{array}
\end{equation}
$$
Cy=\sum_{k=1}^n \sum_{\nu=0}^{k-1} (-1)^{\nu}C_k^\nu \sum_{j=n-k+\nu}^{n-1}y^{(j)}(0)\frac{J^\nu(x^{j-n+k-\nu}P_k)}{(j-n+k-\nu)!} \qquad\qquad\qquad\qquad
$$
$$
\qquad+i\sum_{k=0}^{n-1} \sum_{\nu=0}^k (-1)^{\nu}C_k^\nu
\sum_{j=n-k+\nu-1}^{n-1}y^{(j)}(0)\frac{J^\nu(x^{j-n+k-\nu+1}Q_k)}{(j-n+k-\nu+1)!}
$$
$$
\qquad\qquad\qquad\qquad +i\sum_{k=1}^{n-1} \sum_{\nu=0}^{k-1} (-1)^{\nu}C_k^\nu\sum_{j=n-k+\nu}^{n-1}y^{(j)}(0)\frac{J^{\nu+1}(x^{j-n+k-\nu}Q_k)}{(j-n+k-\nu)!},
$$
откуда, в свою очередь, вытекает (\ref{eq6.3}) и (\ref{eq6.4}). Преобразуя (\ref{eq6.8}) с учетом
\begin{equation}\label{eq6.9}
J^kf=\int_0^x \frac{(x-t)^{k-1}}{(k-1)!}f(t)\,dt, \quad k\in{\mathbb N},
\end{equation}
приходим к представлению для $K(x,t)$ в утверждении теоремы.  $\hfill\Box$
\end{proof}

\begin{remark}\label{rem3}
Если $Q_0(x)=0$ и $P_k(x),Q_k(x)\in C[0,1]$ при $k\ge1,$ то ядро $K(x,t)$ в (\ref{eq6.3}) непрерывно, и $K(x,x)\equiv0.$ Например, это будет иметь место, если на единицу уменьшить сингулярность коэффициентов
выражения (\ref{eq3.10}) и обнулить коэффициент при $(2n-1)$-й производной, т.е. вместо (\ref{eq6.1}) потребовать
$$
p_0, \frac1{p_0}\in L(0,1), \quad p_k\in W_{2,\frac1b}^{1-k}, \;\; k=\overline{1,n}, \quad q_0=0, \quad q_\nu\in W_{2,\frac1b}^{1-\nu}, \;\;\nu=\overline{1,n-1}.
$$
\end{remark}

\begin{remark}\label{rem4}
В то время как функции $P_k$ и $Q_k$ в теореме~\ref{th6} определены с точностью до многочленов степени меньшей~$k,$ ядро $K(x,t)$  не зависит от выбора этих функций. Аналогичное замечание будет справедливым и
применительно к теореме~\ref{th7}, относящейся к выражению нечетного порядка.
\end{remark}

\subsection{Нечетный порядок}
Согласно разделу~\ref{s5.2} условия (\ref{eq3.13}) примут вид
\begin{equation}\label{eq6.10}
\frac1{q_0}\in W_1^1, \quad q_k\in W_1^{1-k}, \;\; p_k\in W_1^{-k}, \;\; k=\overline{0,n},
\end{equation}
что отвечает случаю $p^*=1$ и $b=1.$ Тогда выражение (\ref{eq3.12}) определено в точности на функциях $y\in W_\infty^n$ и является функционалом $\ell_{2n+1}y\in W_\infty^{-n-1}.$

\begin{theorem}\label{th7}
Для выражения (\ref{eq3.12}) с коэффициентами (\ref{eq6.10}) имеет место представление
\begin{equation}\label{eq6.11}
\ell_{2n+1}y=\frac{d^{n+1}}{dx^{n+1}}\Big(By^{(n)}+Cy\Big), \quad 0<x<1,
\end{equation}
где $B$ взаимно однозначно отображает $L_\infty(0,1)$ на $L_\infty(0,1),$ причем
\begin{equation}\label{eq6.12}
Bf=2iq_0(x)f(x) +\int_0^x K(x,t)f(t)\,dt, \quad |K(x,t)|\le a(t)\in L(0,1),
\end{equation}
\begin{equation}\label{eq6.13}
Cy=\sum_{\nu=0}^{n-1} y^{(\nu)}(0)u_\nu(x), \quad u_\nu(x)\in W_1^1[0,1], \quad \nu=\overline{0,n-1}.
\end{equation}
Кроме того, справедливо представление
$$
K(x,t)=\sum_{k=0}^n (-1)^k\frac{(x-t)^k}{k!}(P_k(t)-iQ_k(t))
+\sum_{k=1}^n\int_t^x(P_k(\tau)+iQ_k(\tau))\frac{(\tau-t)^{k-1}}{(k-1)!}\,d\tau
$$
$$
+\sum_{k=2}^n\sum_{\nu=1}^{k-1} (-1)^\nu\int_t^x\frac{(x-\tau)^\nu}{\nu!}\Big(C_k^\nu P_k(\tau)
+i(C_{k-1}^\nu-C_{k-1}^{\nu-1})Q_k(\tau)\Big) \frac{(\tau-t)^{k-1-\nu}}{(k-1-\nu)!}\,d\tau,
$$
где $P_0=p_0$ и $Q_0=q_0',$ а также $P_k^{(k)}=p_k$ и $Q_k^{(k-1)}=q_k$ при $k=\overline{1,n}.$
\end{theorem}

\begin{proof} Перепишем (\ref{eq3.12}) в виде
$$
\ell_{2n+1}y:=i\sum_{k=0}^{n} \Big((Q_k^{(k-1)}y^{(n-k+1)})^{(n-k)} +(Q_k^{(k-1)}y^{(n-k)})^{(n-k+1)}\Big)\qquad\qquad
$$
\begin{equation}\label{eq6.14}
\qquad\qquad\qquad \qquad\qquad\qquad\qquad\qquad\qquad\qquad+ \sum_{k=0}^n (P_k^{(k)}y^{(n-k)})^{(n-k)},
\end{equation}
где $Q_0^{(-1)}=q_0\in W_1^1,$ причем $y^{(n)}\in W_\infty^0.$ Согласно предложению~\ref{pr5} имеем
$$
(Q_0^{(-1)}y^{(n+1)})^{(n)}+(Q_0^{(-1)}y^{(n)})^{(n+1)}= \frac{d^{n+1}}{dx^{n+1}}(2q_0y^{(n)}-J(Q_0y^{(n)})).
$$
Поскольку $y\in W_\infty^n,$ с помощью предложения~\ref{pr7} также получим
$$
Q_k^{(k-1)}y^{(n-k+1)}=\sum_{\nu=0}^{k-1} (-1)^{\nu}C_{k-1}^\nu(Q_k y^{(n-k+\nu+1)})^{(k-\nu-1)}, \quad k=\overline{1,n},
$$
$$
Q_k^{(k-1)}y^{(n-k)}=\sum_{\nu=0}^{k-1} (-1)^{\nu}C_{k-1}^\nu(Q_k y^{(n-k+\nu)})^{(k-\nu-1)}, \quad k=\overline{1,n},
$$
$$
P_k^{(k)}y^{(n-k)}=\sum_{\nu=0}^k (-1)^{\nu}C_k^\nu(P_k y^{(n-k+\nu)})^{(k-\nu)}, \quad k=\overline{0,n}.
$$
Подставляя последние четыре представления в (\ref{eq6.14}), будем иметь
$$
\ell_{2n+1}y=\frac{d^{n+1}}{dx^{n+1}}\Bigg(2i q_0y^{(n)}-iJ(Q_0y^{(n)}) +\sum_{k=0}^n \sum_{\nu=0}^k (-1)^{\nu}C_k^\nu J^{\nu+1}(P_k y^{(n-k+\nu)}) \quad
$$
$$
\qquad\qquad+i\sum_{k=1}^n\sum_{\nu=0}^{k-1} (-1)^{\nu}C_{k-1}^\nu  \Big(J^{\nu+2}(Q_k y^{(n-k+\nu+1)}) +J^{\nu+1}(Q_k y^{(n-k+\nu)})\Big)\Bigg).
$$
При помощи (\ref{eq6.7}) последнее выражение преобразуется к виду (\ref{eq6.11}), где
$$
Bf=2i q_0f -iJ(Q_0f)+\sum_{k=0}^n \sum_{\nu=0}^k (-1)^{\nu}C_k^\nu J^{\nu+1}(P_k J^{k-\nu}f) \qquad\qquad\qquad
$$
\begin{equation}\label{eq6.15}
\qquad+i\sum_{k=1}^n \sum_{\nu=0}^{k-1} (-1)^{\nu} C_{k-1}^\nu \Big(J^{\nu+2}(Q_k J^{k-\nu-1}f) +J^{\nu+1}(Q_k J^{k-\nu}f)\Big),
\end{equation}
$$
Cy=i\sum_{k=2}^n \sum_{\nu=0}^{k-2} (-1)^{\nu}C_{k-1}^\nu \sum_{j=n-k+\nu+1}^{n-1}y^{(j)}(0)\frac{J^{\nu+2}(x^{j-n+k-\nu-1}Q_k)}{(j-n+k-\nu-1)!}\qquad
$$
$$
\qquad\qquad\quad
 +i\sum_{k=1}^n \sum_{\nu=0}^{k-1} (-1)^{\nu}C_{k-1}^\nu\sum_{j=n-k+\nu}^{n-1}y^{(j)}(0)\frac{J^{\nu+1}(x^{j-n+k-\nu}Q_k)}{(j-n+k-\nu)!}
$$
$$
\qquad\qquad\qquad\qquad\qquad +\sum_{k=1}^n \sum_{\nu=0}^{k-1} (-1)^{\nu}C_k^\nu \sum_{j=n-k+\nu}^{n-1}y^{(j)}(0)\frac{J^{\nu+1}(x^{j-n+k-\nu}P_k)}{(j-n+k-\nu)!},
$$
откуда, в свою очередь, вытекает (\ref{eq6.12}) и (\ref{eq6.13}). Наконец, преобразуя (\ref{eq6.15}) с учетом (\ref{eq6.9}), приходим к требуемому представлению для $K(x,t).$  $\hfill\Box$
\end{proof}

\begin{remark}\label{rem5}
Если $P_0(x)=iQ_0(x)$ и $P_k(x),Q_k(x)\in C[0,1]$ при $k=\overline{1,n},$ то ядро $K(x,t)$ в (\ref{eq6.12}) непрерывно и $K(x,x)\equiv0.$ В частности, это имеет место, если является постоянным коэффициент при
старшей производной и нулевым коэффициент при $2n$-й производной, а также на единицу меньшей сингулярность коэффициентов выражения (\ref{eq3.12}), т.е. если вместо (\ref{eq6.10}) потребовать
$$
q_0\equiv\text{const}, \quad p_0\equiv0, \quad q_k\in W_1^{2-k}, \;\; p_k\in W_1^{1-k}, \;\; k=\overline{1,n}.
$$
\end{remark}

\section{Связь локальных квазипроизводных с нелокальными}\label{s7}

\subsection{Суть вопроса и мотивация}\label{s7.1}
Рассмотрим дифференциальное выражение $\ell_Ny$ вида (\ref{eq3.10}) при $N=2n$ с коэффициентами, подчиняющимися (\ref{eq3.11}), и вида (\ref{eq3.12}) при $N=2n+1$ с коэффициентами, отвечающими (\ref{eq3.13}).

Как было установлено в предыдущем параграфе, выражение $\ell_Ny$ можно представить в виде (\ref{eq1.2}). Поэтому желательно иметь явные формулы, выражающие нелокальные квазипроизводные $y^{\langle k\rangle}$
вида (\ref{eq1.7}) через какие-ни\-будь локальные квазипроизводные $y^{[k]}$ и наоборот. В разделе~\ref{s7.3} мы построим некоторый такой набор $y^{[k]},$ выделяя локальные составляющие в $y^{\langle
k\rangle}.$

В разделе~\ref{s7.4} при помощи результатов из~\S~\ref{s4} будет установлено, что всевозможные наборы локальных квазипроизводных для выражения $\ell_Ny$ связаны друг с другом соотношениями (\ref{eq4.5})
посредством соответствующих наборов абсолютно непрерывных функций вида (\ref{eq4.4}), (\ref{eq4.25}). Параллельно будет дано обобщение предложения~\ref{pr1} на случай выражения $\ell_Ny.$ Также будет указан
способ получения матриц ${\bf F}\in{\cal S}_N[0,1],$ являющихся согласованными с одним и тем же выражением $\ell_Ny,$ позволяющий описать все такие матрицы.

Таким образом, получена равноценная альтернатива регуляризации сингулярных дифференциальных выражений $\ell_Ny,$ отправной точкой для которой является преобразование их к операторно-дифферен\-циальному виду
(\ref{eq1.2}).

\subsection{Система интегральных уравнений}\label{s7.2}
 Для дальнейшего удобно дать следующее определение.
Будем говорить, что измеримая матрица-функция
$$
{\bf K}(x,t)=[K_{k,j}(x,t)]_{k,j=0}^{N-1}, \quad 0<t<x\le1,
$$
принадлежит классу ${\cal K}_N,$ если выполнены следующие условия. Все функции $K_{k,j}(x,t)$ определены для каждого $x\in(0,1],$ причем $K_{k,j}(x,\cdot)\in L(0,x).$ Для почти всех $t\in(0,1)$ имеет место
$K_{k,j}(\,\cdot\,,t)\in AC[t,1],$ и выполняются оценки
\begin{equation}\label{eq7.1}
|K_{k,j}(x,t)|\le a(t), \quad \Big|\frac{\partial}{\partial x}K_{k,j}(x,t)\Big|\le a(x)a(t), \quad k,j=\overline{0,N-1},
\end{equation}
где $a\in L(0,1)$ -- своя для каждой ${\bf K}(x,t).$ Наконец, все $K_{k,j}(x,x)\in L(0,1).$

\begin{propos}\label{pr8}
Для любых $\widetilde{\bf K}(x,t),\widehat{\bf K}(x,t)\in{\cal K}_N$ интегральное уравнение
\begin{equation}\label{eq7.2}
{\bf R}(x,t)=\widetilde{\bf K}(x,t)+\int_t^x \widehat{\bf K}(x,\tau){\bf R}(\tau,t)\,d\tau, \quad 0<t<x\le1,
\end{equation}
имеет единственное решение ${\bf R}(x,t)=[R_{k,j}(x,t)]_{k,j=0}^{N-1}\in{\cal K}_N.$
\end{propos}

\begin{proof} Решая уравнение (\ref{eq7.2}) методом
последовательных приближений, будем иметь
\begin{equation}\label{eq7.3}
{\bf R}(x,t)=\sum_{k=0}^\infty {\bf R}_k(x,t),
\end{equation}
где
\begin{equation}\label{eq7.4}
{\bf R}_0(x,t)=\widetilde{\bf K}(x,t), \quad {\bf R}_k(x,t)= \int_t^x \widehat{\bf K}(x,\tau){\bf R}_{k-1}(\tau,t)\,d\tau, \quad k\ge1.
\end{equation}

Пусть одна и та же функция $a\in L(0,1)$ отвечает обеим матричным функциям $\widetilde{\bf K}(x,t)$ и $\widehat{\bf K}(x,t)$ в смысле определения класса ${\cal K}_N$ (см. неравенства (\ref{eq7.1})). Рассмотрим
матричную функцию $A(t)$ той же размерности, всеми элементами которой является $a(t).$ По индукции нетрудно получить матричную оценку
\begin{equation}\label{eq7.5}
|{\bf R}_k(x,t)|\le\frac1{k!}\Big(\int_t^x A(\tau)\,d\tau\Big)^kA(t), \quad k\ge0,
\end{equation}
где элементами матрицы $|{\bf R}_k(x,t)|$ являются модули соответствующих элементов матрицы ${\bf R}_k(x,t),$ а неравенство поэлементно. Следовательно, для каждого фиксированного $x\in(0,1]$ ряд в (\ref{eq7.3})
сходится в $L(0,x)$ равномерно по~$x,$ а значит, его сумма ${\bf R}(x,t)$ является решением уравнения (\ref{eq7.2}).

Далее, в силу соотношения (\ref{eq7.2}) будем иметь $R_{k,j}(\,\cdot\,,t)\in AC[t,1]$ для почти всех~$t\in(0,1).$ Наконец, используя (\ref{eq7.2}), (\ref{eq7.3}) и (\ref{eq7.5}), получаем оценки
$$
|R_{k,j}(x,t)|\le \widetilde a(t), \quad \Big|\frac{\partial}{\partial x}R_{k,j}(x,t)\Big|\le \widetilde a(x)\widetilde a(t), \quad k,j=\overline{0,N-1},
$$
где $\widetilde a\in L(0,1)$ определяется по $a.$ Таким образом, приходим к ${\bf R}(x,t)\in{\cal K}_N.$

Единственность будет следовать из оценки (\ref{eq7.5}), если в уравнении (\ref{eq7.2}) положить $\widetilde{\bf K}(x,t)=0,$ а в (\ref{eq7.4}) заменить $\widetilde{\bf K}(x,t)$ на ${\bf R}(x,t).$ $\hfill\Box$
\end{proof}

\subsection{Согласованная матрица. Пересчет $y^{\langle k\rangle}$ через $y^{[k]}$ и обратно}\label{s7.3}
Следующая теорема дает некоторую согласованную с $\ell_Ny$  матрицу ${\bf F}\in{S}_N[0,1]$ и позволяет явно выразить $y^{\langle k\rangle}$ через порождаемые ею квазипроизводные $y^{[k]}.$

\begin{theorem}\label{th8}
Найдется такая матрица ${\bf F}=[f_{k,j}]_{k,j=1}^N\in{\cal S}_N[0,1],$ согласованная с выражением $\ell_Ny,$ а также матрица ${\bf K}(x,t)=[K_{k,j}(x,t)]_{k,j=0}^{N-1}\in{\cal K}_N,$ что квазипроизводные
(\ref{eq3.9}), порожденные этой матрицей ${\bf F},$ будут связаны с нелокальными квазипроизводными (\ref{eq1.7}) соотношениями
\begin{equation}\label{eq7.6}
y^{\langle k\rangle}(x)=y^{[k]}(x)+\sum_{j=0}^k\int_0^x K_{k,j}(x,t)y^{[j]}(t)\,dt, \quad k=\overline{0,N-1},
\end{equation}
причем $K_{k,j}(x,t)=0$ при $k<\max\{n,j\},$ где $n=[\frac{N}2].$
\end{theorem}

\begin{proof} Естественно положить
$$
y^{[k]}:=y^{(k)}, \quad k=\overline{0,n-1}, \quad f_{k,j}:=\delta_{k,j-1}, \quad k<\max\{n,j-1\}.
$$

1) Пусть $N=2n.$ Определим
\begin{equation}\label{eq7.7}
y^{[n]}:=p_0y^{(n)}+\sum_{j=1}^n (P_j+iQ_{j-1})y^{(n-j)},
\end{equation}
что можно записать в виде (\ref{eq3.9}) при $k=n,$ если положить
\begin{equation}\label{eq7.8}
f_{n,n+1}:=\frac1{p_0}, \quad f_{n,j}:=-\frac{P_{n+1-j}+iQ_{n-j}}{p_0}, \quad j=\overline{1,n}.
\end{equation}
При этом получаем представление
\begin{equation}\label{eq7.9}
y^{(n)}=\sum_{j=1}^{n+1}f_{n,j}y^{[j-1]}.
\end{equation}
В соответствии с формулами (\ref{eq1.7}), (\ref{eq6.2}), (\ref{eq6.6}) и (\ref{eq7.7}) имеем
$$
\begin{array}{c}
\displaystyle y^{\langle n\rangle}=y^{[n]} +\sum_{j=1}^n \sum_{\nu=1}^j (-1)^{\nu}C_j^\nu J^\nu(P_j y^{(n-j+\nu)})
\quad\qquad\qquad\qquad\qquad\\[5mm]
\quad\qquad\qquad\displaystyle + i\sum_{j=1}^{n-1}\sum_{\nu=1}^j (-1)^{\nu}C_j^\nu  J^\nu(Q_j y^{(n-j+\nu-1)})\\[5mm]
\qquad\qquad\qquad\qquad\qquad\qquad\qquad\displaystyle +i\sum_{j=0}^{n-1}\sum_{\nu=0}^j (-1)^{\nu}C_j^\nu  J^{\nu+1}(Q_j
y^{(n-j+\nu)}).
\end{array}
$$
Меняя порядок суммирования, будем иметь
\begin{equation}\label{eq7.10}
y^{\langle n\rangle}=y^{[n]} +\sum_{j=0}^n \int_0^x \widetilde K_{n,j}(x,t) y^{(j)}(t)\,dt,
\end{equation}
где
\begin{equation}\label{eq7.11}
\widetilde K_{n,j}(x,t)=
 \left\{\begin{array}{l}
\quad\;\; 0, \quad j=0,\\[3mm]
\displaystyle \sum_{\nu=n-j+1}^n (-1)^{j-n+\nu}\frac{(x-t)^{j-n+\nu-1}}{(j-n+\nu-1)!}\Big(C_\nu^{j-n+\nu}P_\nu(t)\\[5mm]
\displaystyle \;\;\qquad+i\Big(C_{\nu-1}^{j-n+\nu}-C_{\nu-1}^{j-n+\nu-1}\Big)Q_{\nu-1}(t)\Big), \quad j=\overline{1,n-1},\\[5mm]
\quad \displaystyle\sum_{\nu=1}^n (-1)^\nu \frac{(x-t)^{\nu-1}}{(\nu-1)!}(P_\nu(t)-iQ_{\nu-1}(t)), \quad j=n.
\end{array}\right.
\end{equation}
Подставляя (\ref{eq7.9}) в (\ref{eq7.10}), получаем (\ref{eq7.6}) для $k=n,$ т.е.
\begin{equation}\label{eq7.12}
y^{\langle n\rangle}=y^{[n]} +\sum_{j=0}^n \int_0^x K_{n,j}(x,t) y^{[j]}(t)\,dt,
\end{equation}
где
\begin{equation}\label{eq7.13}
K_{n,j}(x,t)=(1-\delta_{n,j})\widetilde K_{n,j}(x,t) +\widetilde K_{n,n}(x,t) f_{n,j+1}(t), \quad j=\overline{0,n}.
\end{equation}

Если $n>1,$ то предположим по индукции, что при $k\le s-1$ для некоторого $s\in\{n+1,\ldots,2n-1\}$ построены все $K_{k,j}(x,t)$ и $f_{k,j},$ а значит, при таких~$k$ определены и $y^{[k]}.$ Пусть, кроме того,
имеют место соотношения
\begin{equation}\label{eq7.14}
y^{\langle s-1\rangle}=y^{[s-1]} +\sum_{j=0}^n \int_0^x K_{s-1,j}(x,t) y^{[j]}(t)\,dt,
\end{equation}
\begin{equation}\label{eq7.15}
K_{s-1,j}(x,t)=(1-\delta_{n,j})\widetilde K_{s-1,j}(x,t) +\widetilde K_{s-1,n}(x,t) f_{n,j+1}(t), \quad j=\overline{0,n},
\end{equation}
\begin{equation}\label{eq7.16}
\widetilde K_{s-1,j}(x,t)= \left\{\begin{array}{l}
\quad\;\; 0, \quad j=\overline{0,s-n-1},\\[2mm]
\displaystyle \sum_{\nu=s-j}^n (-1)^{j-n+\nu}\frac{(x-t)^{j-s+\nu}}{(j-s+\nu)!}\Big(C_\nu^{j-n+\nu}P_\nu(t)\\[5mm]
\displaystyle \;\;+i\Big(C_{\nu-1}^{j-n+\nu}-C_{\nu-1}^{j-n+\nu-1}\Big)Q_{\nu-1}(t)\Big), \; j=\overline{s-n,n-1},\\[5mm]
\displaystyle\sum_{\nu=s-n}^n (-1)^\nu \frac{(x-t)^{n-s+\nu}}{(n-s+\nu)!}(P_\nu(t)-iQ_{\nu-1}(t)), \quad j=n,
\end{array}\right.
\end{equation}
совпадающие при $s=n+1$ с (\ref{eq7.12}), (\ref{eq7.13}) и (\ref{eq7.11}), соответственно.

Дифференцируя (\ref{eq7.14}) и учитывая (\ref{eq1.7}), получаем
\begin{equation}\label{eq7.17}
y^{\langle s\rangle}=y^{[s]} +\sum_{j=0}^n \int_0^x K_{s,j}(x,t) y^{[j]}(t)\,dt,
\end{equation}
где $y^{[s]}$ определяется формулой (\ref{eq3.9}) при $k=s,$ в которой
$$
f_{s,j}:=-K_{s-1,j-1}(x,x), \quad j=\overline{1,n+1}, \quad f_{s,j}:=\delta_{s+1,j}, \quad j>n+1;
$$
$$
K_{s,j}(x,t)=\frac{\partial}{\partial x} K_{s-1,j}(x,t), \quad j=\overline{0,n}, \quad K_{s,j}(x,t):=0, \quad j>n.
$$
Используя (\ref{eq7.15}) и (\ref{eq7.16}), вычисляем
$$
f_{s,j}=(\delta_{n,j}-1)\widetilde f_{s,j-1} -\widetilde f_{s,n} f_{n,j}, \quad j=\overline{1,n+1},
$$
где
$$
\widetilde f_{s,j}=
 \left\{\begin{array}{l}
\quad\;\; 0, \quad j=\overline{0,s-n-1},\\[3mm]
\displaystyle (-1)^{s-n}\Big(C_{s-j}^{s-n}P_{s-j}+i\Big(C_{s-j-1}^{s-n}-C_{s-j-1}^{s-n-1}\Big)Q_{s-j-1}\Big), \;\;  j=\overline{s-n,n-1},\\[4mm]
\displaystyle (-1)^{s-n} (P_{s-n}-iQ_{s-n-1}), \quad j=n,
\end{array}\right.
$$
а также
\begin{equation}\label{eq7.18}
K_{s,j}(x,t)=(1-\delta_{n,j})\widetilde K_{s,j}(x,t) +\widetilde K_{s,n}(x,t) f_{n,j+1}(t), \quad j=\overline{0,n},
\end{equation}
где
\begin{equation}\label{eq7.19}
\widetilde K_{s,j}(x,t)=
 \left\{\begin{array}{l}
\quad\;\; 0, \quad j=\overline{0,s-n},\\[3mm]
\displaystyle \sum_{\nu=s-j+1}^n (-1)^{j-n+\nu}\frac{(x-t)^{j-s+\nu-1}}{(j-s+\nu-1)!}\Big(C_\nu^{j-n+\nu}P_\nu(t)\\[5mm]
\displaystyle +i\Big(C_{\nu-1}^{j-n+\nu}-C_{\nu-1}^{j-n+\nu-1}\Big)Q_{\nu-1}(t)\Big), \; j=\overline{s-n+1,n-1},\\[4mm]
\displaystyle\sum_{\nu=s-n+1}^n (-1)^\nu \frac{(x-t)^{n-s+\nu-1}}{(n-s+\nu-1)!}(P_\nu(t)-iQ_{\nu-1}(t)), \;\; j=n.
\end{array}\right.
\end{equation}
Поскольку формулы (\ref{eq7.17}), (\ref{eq7.18}) и (\ref{eq7.19}) могут быть получены соответственно из (\ref{eq7.14}), (\ref{eq7.15}) и (\ref{eq7.16}) заменой $s$ на $s+1,$ все предположения индукции являются
верными и в силу произвольности $s$ матрица ${\bf K}(x,t)$ построена.

Остается найти последнюю строку матрицы ${\bf F}.$ При $s=2n-1$ соотношения (\ref{eq7.18}) и (\ref{eq7.19}) принимают вид
$$
K_{2n-1,j}(x,t)=(1-\delta_{n,j})\widetilde K_{2n-1,j}(x,t) +\widetilde K_{2n-1,n}(x,t) f_{n,j+1}(t), \quad j=\overline{0,n},
$$
$$
\widetilde K_{2n-1,j}(x,t)=
 \left\{\begin{array}{l}
\quad\;\; 0, \quad j=\overline{0,n-1},\\[3mm]
\displaystyle (-1)^n (P_n(t)-iQ_{n-1}(t)), \quad j=n.
\end{array}\right.
$$
Следовательно, дифференцируя (\ref{eq7.17}) при $s=2n-1,$ с учетом (\ref{eq1.7}) и~(\ref{eq6.2}) будем иметь (\ref{eq3.8}) при $N=2n,$ где ${\cal F}y=\ell_{2n}y$ и
\begin{equation}\label{eq7.20}
f_{2n,j}=(-1)^n(iQ_{n-1}-P_n)f_{n,j}, \quad j=\overline{1,n+1}, \quad f_{2n,j}=0, \quad j>n+1.
\end{equation}
Построенная матрица ${\bf F}$ в силу (\ref{eq3.11}) принадлежит классу ${\cal S}_{2n}[0,1]$ и порождает квазидифференциальное выражение ${\cal F}y,$ совпадающее с $\ell_{2n}y$ для всех $y\in{\cal D}({\cal F}),$
что и завершает доказательство теоремы для случая $N=2n.$

2) Случай $N=2n+1$ рассматривается аналогично. При этом формулы для нахождения элементов матриц ${\bf F}$ и ${\bf K}(x,t)$ будут иметь вид
$$
f_{n,j}=\frac{\delta_{n+1,j}}{2iq_0}, \quad j=\overline{1,2n+1},
$$
$$
f_{k,j}=
 \left\{\begin{array}{l}
\quad\;\; 0, \quad j=\overline{1,k-n-1},\\[3mm]
\displaystyle (-1)^{k-n}\Big(C_{k-j}^{k-n-1}P_{k-j}
+i\Big(C_{k-j-1}^{k-n-1}-C_{k-j-1}^{k-n-2}\Big)Q_{k-j}\Big), \quad  j=\overline{k-n,n},\\[3mm]
\displaystyle (-1)^{k-n} \frac{P_{k-n-1}-iQ_{k-n-1}}{2iq_0}, \quad j=n+1,\\[3mm]
\quad\delta_{k+1,j}, \quad j=\overline{n+2,2n+1},
\end{array}\right.
$$
где $k=\overline{n+1,2n+1},$ а также введено обозначение $C_\nu^{-1}:=0$ при $\nu\ge0;$ и
$$
K_{k,j}(x,t)=
 \left\{\begin{array}{l}
\quad\;\; 0, \quad j\ne\overline{k-n,n},\\[3mm]
\displaystyle \sum_{\nu=k-j}^n (-1)^{j-n+\nu}\frac{(x-t)^{j-k+\nu}}{(j-k+\nu)!}\Big(C_\nu^{j-n+\nu}P_\nu(t)\\[3mm]
\displaystyle \qquad\qquad +i\Big(C_{\nu-1}^{j-n+\nu}-C_{\nu-1}^{j-n+\nu-1}\Big)Q_\nu(t)\Big), \quad j=\overline{k-n,n-1},\\[3mm]
\displaystyle \frac1{2iq_0(t)}\sum_{\nu=k-n}^n (-1)^\nu \frac{(x-t)^{n-k+\nu}}{(n-k+\nu)!}(P_\nu(t)-iQ_\nu(t)), \quad j=n,
\end{array}\right.
$$
где $k=\overline{n,2n}.$ При $k<n$ функции $f_{k,j}$ и $K_{k,j}(x,t)$ были определены выше. $\hfill\Box$
\end{proof}

{\sc Следствие 2. }{\it Представления (\ref{eq7.6}) равносильны представлениям
\begin{equation}\label{eq7.21}
y^{[k]}(x)=y^{\langle k\rangle}(x)+\sum_{j=0}^k\int_0^x R_{k,j}(x,t)y^{\langle j\rangle}(t)\,dt, \quad k=\overline{0,N-1},
\end{equation}
причем ${\bf R}(x,t):=[R_{k,j}(x,t)]_{k,j=0}^{N-1}\in{\cal K}_N,$ где $R_{k,j}(x,t)=0$ при $k<\max\{n,j\}.$}

\medskip
{\sc Доказательство.} В силу предложения~\ref{pr8}, для ядра ${\bf K}(x,t)\in{\cal K}_N$ в системе уравнений (\ref{eq7.6}) относительно $y^{[k]}$ существует резольвентное ядро ${\bf R}(x,t)\in{\cal K}_N,$
удовлетворяющее уравнению (\ref{eq7.2}) при $\widetilde{\bf K}(x,t)=\widehat{\bf K}(x,t)=-{\bf K}(x,t)$ и являющееся нижнетреугольным, как и исходное. $\hfill\Box$

\begin{remark}\label{rem6}
Доказательство теоремы~\ref{th8} конструктивно и дает явные формулы для нахождения согласованной с выражением $\ell_Ny$ матрицы ${\bf F}\in{S}_N[0,1],$ а также ядра ${\bf K}(x,t)$ интегрального преобразования
(\ref{eq7.6}).
\end{remark}

{\sc Пример 1.} Пусть $\ell_2y=-ly,$ где $ly$ имеет вид (\ref{eq3.1}). Тогда в соответствии с (\ref{eq3.2}) и  (\ref{eq3.10}) имеем $n=1,$ $p_0=1,$ $P_1=-\sigma$ и $Q_0=0.$ Теорема~\ref{th6} дает $u_0=-\sigma$
и $K(x,t)=\sigma(t)-\sigma(x),$ причем $b=1.$ Согласно (\ref{eq1.7}), (\ref{eq6.3}) и (\ref{eq6.4}) вычисляем
$$
y^{\langle 1\rangle}=y'(x)+\int_0^x(\sigma(t)-\sigma(x))y'(t)\,dt-y(0)\sigma(x)=y^{[1]}(x)+\int_0^x\sigma(t)y'(t)\,dt,
$$
где квазипроизводная $y^{[1]}$ определяется формулой (\ref{eq3.4}). Кроме того, согласно (\ref{eq7.8}) и (\ref{eq7.20}) соответствующая матрица ${\bf F} =[f_{k,j}]_{k,j=1}^{2}$ совпадает с ${\bf Q}$ в
(\ref{eq3.6}).

\subsection{Об общем виде согласованных матриц}\label{s7.4}
Теоремы~\ref{th4} и~\ref{th5} позволяют найти и другие матрицы $\widetilde{\bf F}\in {\cal S}_N[0,1],$ согласованные с одним и тем же выражением $\ell_N y$ по построенной в предыдущем разделе матрице ${\bf F}.$
Для этого нужно задаться всевозможными наборами абсолютно непрерывных функций (\ref{eq4.4}), (\ref{eq4.25}). При этом элементы соответствующих $\widetilde{\bf F}$ определяются по формулам (\ref{eq4.13}) c
$g_{0,0}\equiv1$ и из треугольных алгебраических систем (\ref{eq4.14}) и (\ref{eq4.26}).

Кроме того, теперь мы можем доказать следующее обобщение предложения~\ref{pr1} на случай произвольного выражения $\ell_Ny,$ из
которого также следует, что указанным способом определяются {\bf все} матрицы, согласованные с $\ell_N y.$

\begin{theorem}\label{th9}
Все наборы локальных квазипроизводных $y^{\{k\}},$ $k=\overline{1,N-1},$ для выражения $\ell_Ny$ вида (\ref{eq3.10}) с четным $N$ и вида (\ref{eq3.12}) -- с нечетным образуют $\pi_N$-па\-ра\-ме\-трическое
семейство (\ref{eq4.5}), подчиненное $\pi_N:=(1+N)N/2-2$ функциональным параметрам (\ref{eq4.4}) (при том, что выполняется (\ref{eq4.25})).
\end{theorem}

\begin{proof}
В разделе~\ref{s7.3} построена согласованная с выражением $\ell_Ny$ матрица ${\bf F}\in{S}_N[0,1].$ Очевидно, все эквивалентные ей матрицы также являются согласованными с $\ell_Ny,$ а в силу теорем~\ref{th4}
и~\ref{th5} порождаемые ими наборы квазипроизводных образуют указанное $\pi_N$-па\-ра\-метрическое семейство (\ref{eq4.5}).

Таким образом, остается показать, что у $\ell_Ny$ нет иных согласованных матриц. Пусть $\widetilde{\bf F}\in{S}_N[0,1]$ является согласованной с $\ell_Ny.$ Тогда для всех $y\in{\cal D}(\widetilde{\cal F})$
имеем $\ell_Ny\in L(0,1).$ Согласно теоремам~\ref{th6} и~\ref{th7}, а также~(\ref{eq1.7}) последнее влечет $y^{\langle k\rangle}\in AC[0,1]$ при $k=\overline{0,N-1}.$ В силу следствия~2 в предыдущем разделе
получаем $y\in{\cal D}({\cal F}),$ т.е. ${\cal D}(\widetilde{\cal F})\subset{\cal D}({\cal F}).$ Тогда следствие~1 в разделе~\ref{s4.2} дает ${\cal D}(\widetilde{\cal F})={\cal D}({\cal F}),$ т.е. матрицы ${\bf
F}$ и $\widetilde{\bf F}$ являются полуэквивалентными, а значит, и эквивалентными, поскольку ${\cal F}y=\ell_Ny=\widetilde{\cal F}y$ для всех $y\in{\cal D}({\cal F}).$ $\hfill\Box$
\end{proof}

{\sc Пример 2.} Покажем, что все матрицы $[\widetilde f_{k,j}]_{k,j=1}^2,$ согласованные с выражением $-ly$ вида (\ref{eq3.1}), действительно имеют вид (\ref{eq3.7}). Для этого роль модельной матрицы ${\bf F}$
отведем ${\bf Q}$ из (\ref{eq3.6}). При $N=2$ функции (\ref{eq4.4}), (\ref{eq4.25}) имеют вид $g_{1,0}\in AC[0,1]$ и $g_{1,1}\equiv1.$ Тогда (\ref{eq4.13}) дает $\widetilde f_{1,2}=1.$ Система (\ref{eq4.14})
приводит непосредственно к формуле $\widetilde f_{1,1}=-g_{1,0}+\sigma,$ а система (\ref{eq4.26}) принимает вид
$$
\left[\begin{array}{cc}1& g_{1,0}\\[3mm]
0 & 1
\end{array}\right]
\left[\begin{array}{c} \widetilde f_{2,1}\\[3mm]
\widetilde f_{2,2}
\end{array}\right]=
\left[\begin{array}{c} g_{1,0}'+g_{1,0}\sigma-\sigma^2\\[3mm]
g_{1,0}-\sigma
\end{array}\right],
$$
откуда получаем $\widetilde f_{2,2}=g_{1,0}-\sigma$ и
$$
\widetilde f_{2,1}=g_{1,0}'+g_{1,0}\sigma-\sigma^2 -g_{1,0}\widetilde f_{2,2} =g_{1,0}'-(\sigma-g_{1,0})^2.
$$
Таким образом, для всякой функции $g_{1,0}\in AC[0,1]$ приходим к согласованной с выражением $-ly$ матрице (\ref{eq3.7}), в которой $\sigma_1=\sigma-g_{1,0}$ и $q_1=g_{1,0}'.$

\medskip
Возможность выражать явно любые наборы квазипроизводных, включая нелокальные, друг через друга потребуется в~\S~\ref{s9} при рассмотрении краевой задачи для сингулярного дифференциального уравнения
$\ell_Ny=\lambda y.$

\section{Доказательство теоремы~\ref{th1}}\label{s8}

\subsection{Вводное замечание}
В настоящем параграфе строится оператор $A,$ являющийся обратным к оператору $L,$ порожденному выражением (\ref{eq1.2}) и краевыми условиями (\ref{eq1.5}), (\ref{eq1.10}) при $l>N-l>0,$ где $N=m+n.$

При этом предполагается, что операторы $B$ и $C$ имеют вид (\ref{eq1.4}) и (\ref{eq1.3}), соответственно, где функция $K(x,t)$ непрерывна в треугольнике $0\le t\le x\le 1$ и удовлетворяет условию
$K(x,x)\equiv0,$ тогда как $u_\nu(x)\in L(0,1),$ $\nu=\overline{0,n-1}.$

Мы покажем, что построенный оператор $A$ имеет вид (\ref{eq2.3}), и для него будут выполнены условия теоремы~\ref{th3}. Тогда, поскольку система с.п.ф. оператора $L$ является таковой и для  $A,$ утверждение
теоремы~\ref{th1} вытекает из теоремы~\ref{th3}.

Заметим, что не ограничивает общности предположение об обратимости оператора~$L.$ В самом деле, выражение $\widetilde\ell y:=\ell y -\lambda y$ тоже имеет вид~(\ref{eq1.2})--(\ref{eq1.4}):
$$
\widetilde\ell y=\frac{d^m}{dx^m}\Big((I+\widetilde K)y^{(n)}+\widetilde Cy\Big), \;\, \widetilde K:=K-\lambda J^N, \;\, \widetilde Cy:=Cy-\lambda \sum_{j=0}^{n-1}y^{(j)}(0)\frac{x^{j+m}}{(j+m)!},
$$
а порождаемый им оператор с теми же краевыми условиями (\ref{eq1.5}), (\ref{eq1.10}) обладает той же системой с.п.ф., что и $L.$ При этом данные условия в терминах квазипроизводных, порождаемых  по аналогии с
(\ref{eq1.7}) выражением $\widetilde\ell y,$ будут иметь тот же вид. Кроме того, можно показать, что спектр оператора $L$ является дискретным, а значит, найдется $\lambda$ такое, что оператор $L-\lambda I$
обратим.

\subsection{Общее решение уравнения $\ell y=f$}
Операторы
\begin{equation}\label{eq8.1}
R:=(I+K)^{-1}-I, \quad M:=J^n(I+R)J^m
\end{equation}
имеют вид
$$
Rf=\int_0^x R(x,t)f(t)\,dt, \quad Mf=\int_0^x M(x,t)f(t)\,dt,
$$
где ядро $R(x,t)$ непрерывно и является решением интегрального уравнения
\begin{equation}\label{eq8.2}
K(x,t)+R(x,t)+\int_t^x K(x,\tau)R(\tau,t)\,d\tau=0,
\end{equation}
в то время как для $M(x,t)$ справедливо представление
\begin{equation}\label{eq8.3}
M(x,t)=\frac{(x-t)^{N-1}}{(N-1)!}+\int_t^x \frac{(x-\tau)^{n-1}}{(n-1)!} \,d\tau \int_t^\tau R(\tau,\xi)\frac{(\xi-t)^{m-1}}{(m-1)!}\,d\xi.
\end{equation}

Также введем обозначения
$$
\xi_j(x):=(I+R)\frac{x^j}{j!}, \quad y_{n+j+1}(x):=J^n\xi_j(x), \quad j=\overline{0,m-1},
$$
$$
\eta_\nu(x):=-(I+R)u_\nu(x), \quad y_{\nu+1}(x):=\frac{x^\nu}{\nu!}+J^n\eta_\nu(x), \quad  \nu=\overline{0,n-1}.
$$

\begin{propos}\label{pr10}
Функции $y_j(x),$ $j=\overline{1,N},$ образуют фундаментальную систему решений уравнения $\ell y=0,$ а его общее решение имеет вид
\begin{equation}\label{eq8.4}
y(x)=\sum_{j=1}^Ny^{\langle j-1\rangle}(0)y_j(x).
\end{equation}
Кроме того, справедливы асимптотики
\begin{equation}\label{eq8.5}
y_j(x)=\frac{x^{j-1}}{(j-1)!}(1+o(1)), \quad x\to0, \quad j=\overline{1,N}.
\end{equation}
\end{propos}

\begin{proof} Согласно (\ref{eq1.7}) уравнение $\ell y=0$ принимает вид
$$
(I+K)y^{(n)}+Cy=\sum_{j=0}^{m-1}y^{\langle n+j\rangle}(0)\frac{x^j}{j!}.
$$
В силу (\ref{eq1.3}) и первого соотношения в (\ref{eq8.1}) будем иметь
$$
y^{(n)}=\sum_{j=0}^{m-1}y^{\langle n+j\rangle}(0)(I+R)\frac{x^j}{j!}-\sum_{\nu=0}^{n-1} y^{(\nu)}(0)(I+R)u_\nu.
$$
Интегрируя последнее равенство $n$ раз, получаем
$$
J^ny^{(n)}=y-\sum_{\nu=0}^{n-1} y^{(\nu)}(0)\frac{x^\nu}{\nu!}=\sum_{j=0}^{m-1}y^{\langle n+j\rangle}(0)J^n \xi_j+\sum_{\nu=0}^{n-1} y^{(\nu)}(0)J^n\eta_\nu,
$$
что равносильно (\ref{eq8.4}). Далее, согласно (\ref{eq6.9}) и определению $y_j(x)$ получаем
$$
\Big|\frac{\nu!}{x^\nu}y_{\nu+1}(x)-1\Big|\le \frac{\nu!}{x^\nu}\int_0^x\frac{(x-t)^{n-1}}{(n-1)!}|\eta_\nu(t)|\,dt
\le\int_0^x|\eta_\nu(t)|\,dt, \quad \nu=\overline{0,n-1},
$$
тогда как при $j=\overline{0,m-1}$ будем иметь
$$
y_{j+n+1}(x)=J^n\frac{x^j}{j!}+J^nR\frac{x^j}{j!}, \quad J^n\frac{x^j}{j!}=\frac{x^{j+n}}{(j+n)!},
$$
$$
\Big|\frac{(j+n)!}{x^{j+n}}y_{j+n+1}(x)-1\Big| \qquad\qquad\qquad\qquad\qquad\qquad\qquad\qquad\qquad\qquad\qquad\quad\;
$$
$$
\qquad\qquad\qquad\le \frac{(j+n)!}{x^{j+n}}\int_0^x\frac{(x-t)^{n-1}}{(n-1)!}dt\int_0^t|R(t,\tau)|\frac{\tau^j}{j!}\,d\tau =o(x), \quad
x\to0,
$$
поскольку в силу (\ref{eq8.2}) имеем $R(x,x)\equiv0.$ Таким образом, приходим к (\ref{eq8.5}). $\hfill\Box$
\end{proof}

\begin{remark}\label{rem7}
Итак, общее решение уравнения $\ell y=f\in L(0,1)$ имеет вид
$$
y=Mf+\sum_{j=1}^N C_j y_j.
$$
Причем $y,$ в силу (\ref{eq8.4}), является решением задачи Коши с начальными условиями $y^{\langle j\rangle}(0)=C_{j+1},$ $j=\overline{0,N-1},$ поскольку значения $y^{\langle j\rangle}(0)$ совпадают
с~$y_0^{\langle j\rangle}(0).$ Последнее, в свою очередь, видно из представлений
\begin{equation}\label{eq8.6}
(Mf)^{\langle j\rangle}=(Mf)^{(j)}=\int_0^x\frac{\partial^j}{\partial x^j}M(x,t)f(t)\,dt, \quad j=\overline{0,n-1},
\end{equation}
\begin{equation}\label{eq8.7}
(Mf)^{\langle n+s\rangle}=\frac{d^s}{dx^s}\Big((I+K)(Mf)^{(n)}+CMf\Big)=J^{m-s}f, \quad s=\overline{0,m-1}.
\end{equation}
\end{remark}

\subsection{Построение оператора $A$}
Положим $d:=N-l.$ Пусть $g_k(x),$ $k=\overline{1,d},$ -- линейно независимые решения уравнения $\ell y=0,$ удовлетворяющие краевым условиям (\ref{eq1.5}). Можно также считать, что выполняются условия
\begin{equation}\label{eq8.8}
U_{l+j}(g_k)=-\delta_{j,k}, \quad j,k=\overline{1,d}.
\end{equation}
В самом деле, поскольку нуль не является собственным значением оператора $L,$ имеем $\det[U_{l+j}(g_k)]_{j,k=1}^d\ne0.$ Поэтому в линейной оболочке функций $g_k$ найдется базис $\widetilde g_k,$ для которого
$-[U_{l+j}(\widetilde g_k)]_{j,k=1}^d$ -- единичная матрица.

Согласно замечанию~\ref{rem7} и в силу условия (\ref{eq8.8}) оператор $A=L^{-1}$ имеет вид
\begin{equation}\label{eq8.9}
Af=Mf+ \sum_{k=1}^d U_{l+k}(Mf)g_k.
\end{equation}

Поскольку $R(x,t)$ непрерывна и $R(x,x)\equiv0,$ формула (\ref{eq8.3}) дает оценки
\begin{equation}\label{eq8.10}
M_\nu(x,t):=\frac{\partial^\nu}{\partial x^\nu}M(x,t)=\frac{(x-t)^{N-1-\nu}}{(N-1-\nu)!}(1+o(x-t)), \quad t\to x-0,
\end{equation}
выполняющиеся для $\nu=\overline{0,n-1}$ и равномерно по $x\in(0,1].$

Подставляя $y=Mf$ в условия (\ref{eq1.10}) с учетом (\ref{eq8.6}) и (\ref{eq8.7}) и меняя порядок интегрирования будем иметь
\begin{equation}\label{eq8.11}
U_{l+k}(Mf)=\int_0^1 f(t)v_k(t)\,dt, \quad k=\overline{1,d},
\end{equation}
где
$$
v_k(t)=\left\{\begin{array}{cc}
\displaystyle M_{\sigma_{l+k}}(1,t)+\int_t^1 M_{\sigma_{l+k}}(\tau,t)\,d\Phi_{l+k}(\tau), & \sigma_{l+k}<n,\\[5mm]
\displaystyle \frac{(1-t)^{N-1-\sigma_{l+k}}}{(N-1-\sigma_{l+k})!} +\int_t^1\frac{(\tau-t)^{N-1-\sigma_{l+k}}}{(N-1-\sigma_{l+k})!}\,d\Phi_{l+k}(\tau), & \sigma_{l+k}\ge n.
\end{array}\right.
$$

В силу теоремы о среднем для интеграла Стилтьеса (см., например, \cite[Гл.~VI, \S~6, п.~4]{KolmFom}) для всякой непрерывной функции $g(\tau)$ выполняется оценка
$$
\Big|\int_t^1g(\tau)\,d\Phi_\nu(\tau)\Big|\le\max_{t\le\tau\le1}|g(\tau)|V_t^1[\Phi_\nu],
$$
где $V_t^1[\Phi_\nu]$ -- полная вариация $\Phi_\nu$ на отрезке $[t,1].$
 Поскольку рассматриваемые функции $\Phi_\nu(\tau)$ непрерывны в точке~$\tau=1,$ имеем $V_t^1[\Phi_\nu]\to 0$ при $t\to1.$

Таким образом, согласно (\ref{eq8.10}) приходим к асимптотике
\begin{equation}\label{eq8.12}
v_k(t)=\frac{(1-t)^{N-1-\sigma_{l+k}}}{(N-1-\sigma_{l+k})!}(1+o(1)), \quad t\to1, \quad k=\overline{1,d}.
\end{equation}

Итак, учитывая (\ref{eq8.9}) и (\ref{eq8.11}), мы показали, что оператор $A$ имеет вид (\ref{eq2.3}), причем выполняются условия 1) и 2) теоремы~\ref{th3}. Согласно (\ref{eq8.12}) вместе с (\ref{eq1.9}),
функции $v_k(t),$ $k=\overline{1,d},$ удовлетворяют ее условию~3).

Наконец, в соответствии с предложением~\ref{pr10} существуют такие натуральные числа $b_k\le N,$ $k=\overline{1,d},$ что имеют место представления
$$
g_k(x)=\sum_{\nu=b_k}^N \gamma_{k,\nu} y_\nu(x)=\frac{x^{b_k-1}}{(b_k-1)!}(\gamma_{k,b_k}+o(1)), \quad x\to0, \quad \gamma_{k,b_k}\ne0, \quad k=\overline{1,d}.
$$
Покажем, что величины $b_k$ можно считать разными для разных $k.$ Для этого воспользуемся приемом из \cite{Khrom-04}. Пусть $b_i=b_j$ для некоторых $i<j.$ Тогда
$$
\sum_{k=1}^d g_k(x)v_k(t)=\widetilde g_i(x)v_i(t)+g_j(x)\widetilde v_j(t)+\sum_{{k\ne i,j}\atop{k=1}}^d g_k(x)v_k(t),
$$
где
$$
\widetilde g_i(x)=g_i(x)-\alpha g_j(x), \quad \widetilde v_j(t)=v_j(t)+\alpha v_i(t), \quad \alpha:=\frac{\gamma_{i,b_i}}{\gamma_{j,b_j}},
$$
причем в силу (\ref{eq1.9}), (\ref{eq8.12}) будем иметь
$$
\widetilde v_j(t)\sim v_j(t), \quad t\to1.
$$
Поскольку $\widetilde g_i(x)$ тоже является решением уравнения $\ell y=0,$ то
$$
\widetilde g_i(x)=\sum_{\nu=\widetilde b_i}^N \widetilde\gamma_{i,\nu} y_\nu(x)=\frac{x^{\widetilde b_i-1}}{(\widetilde b_i-1)!}(\widetilde\gamma_{i,\widetilde b_i}+o(1)), \quad x\to0, \quad
\widetilde\gamma_{i,\widetilde b_i}\ne0,
$$
причем $\widetilde b_i>b_j.$ Переопределяя $b_i:=\widetilde b_i,$ $g_i(x):=\widetilde g_i(x)$ и применяя эту процедуру конечное число раз, получаем $b_d<\ldots<b_1.$ Поэтому без ущерба для общности можно
считать, что и $g_k(x),$ $k=\overline{1,d},$ удовлетворяют условию 3) теоремы~\ref{th3}.

\section{Теорема о полноте для сингулярных дифференциальных операторов с нерегулярными полураспадающимися краевыми условиями}\label{s9}

Для краткости введем классы $C^{-k}[0,1]$ распределений $f,$ у которых $k$-я первообразная эквивалентна непрерывной функции. Другими словами,
$$
(f,\varphi)=(-1)^k\int_0^1\sigma(x)\varphi^{(k)}(x)\,dx, \quad \forall\varphi\in C_0^\infty(0,1),
$$
с некоторой функцией $\sigma\in C[0,1].$

Пусть выражение $\ell_Ny$ при $N=2n$ имеет вид~(\ref{eq3.10}) с коэффициентами
\begin{equation}\label{eq9.1}
p_0=1, \;\; p_k\in C^{-k}[0,1], \;\; k=\overline{1,n}, \quad q_0=0, \;\; q_k\in C^{-k}[0,1], \;\; \nu=\overline{1,n-1},
\end{equation}
а при $N=2n+1$ -- вид (\ref{eq3.12}), в котором
\begin{equation}\label{eq9.2}
q_0=\frac1{2i}, \quad p_0=0, \quad q_k\in C^{1-k}[0,1], \;\; p_k\in C^{-k}[0,1], \quad k=\overline{1,n}.
\end{equation}

Рассмотрим краевую задачу для уравнения
\begin{equation}\label{eq9.3}
\ell_Ny(x)=\lambda y(x), \quad 0<x<1,
\end{equation}
с линейно независимыми полураспадающимися краевыми условиями
\begin{equation}\label{eq9.4}
\sum_{k=0}^{N-1}\alpha_{j,k}y^{[k]}(0)=0, \quad j=\overline{1,l},
\end{equation}
\begin{equation}\label{eq9.5}
\sum_{k=0}^{N-1}\Big(\alpha_{j,k}y^{[k]}(0)+\beta_{j,k}y^{[k]}(1)\Big)=0, \quad j=\overline{l+1,N},
\end{equation}
где $l>N-l>0,$ а $y^{[k]}$ образуют какой-либо набор квазипроизводных для выражения $\ell_Ny,$ причем матрица $[\beta_{j,k}]$ имеет максимальный ранг, т.е.
\begin{equation}\label{eq9.6}
{\rm rank}[\beta_{j,k}]=N-l \quad (j=\overline{l+1,N}, \; k=\overline{0,N-1}).
\end{equation}

Например, при $N=2n$ в качестве $y^{[k]}$ можно взять квазипроизводные, порожденные согласованной матрицей из \cite{MirShk-16}, а при $N=2n+1$ -- согласованной матрицей из \cite{MirShk19-arXiv}. Вообще,
подойдет любая матрица, построенная в соответствии с процедурой, описанной в разделе~\ref{s7.4} и позволяющей получить все матрицы, согласованные с выражением $\ell_Ny.$ При этом в силу теоремы~\ref{th9}
краевые условия всякий раз можно преобразовать при помощи соотношений (\ref{eq4.5}).

Также вместо $y^{[k]}$ можно взять нелокальные квазипроизводные $y^{\langle k\rangle},$ если в соответствии с теоремами~\ref{th6} и~\ref{th7} представить $\ell_Ny$ в виде~(\ref{eq1.2}). В любом случае условия
(\ref{eq9.4}) и (\ref{eq9.5}) заданы корректно, поскольку соотношение (\ref{eq9.3}) влечет абсолютную непрерывность участвующих в них квазипроизводных.

Следующее утверждение вытекает из теорем~\ref{th1},~\ref{th6},~\ref{th7},~\ref{th9} и следствия~2 в \S~\ref{s7}.

\begin{theorem}\label{th10}
Система с.п.ф. краевой задачи (\ref{eq9.3})--(\ref{eq9.6}) полна в $L_2(0,1).$
\end{theorem}

\begin{proof} В силу теорем~\ref{th6} и~\ref{th7} для выражения $\ell_Ny$ при условиях на коэффициенты (\ref{eq9.1}) -- для четных $N$ -- и (\ref{eq9.2}) -- для нечетных -- имеет место
представление (\ref{eq1.2}), в котором операторы $B$ и $C$ удовлетворяют условиям теоремы~\ref{th1}. Кроме того, согласно теореме~\ref{th9}, без ущерба для общности можно считать, что квазипроизводные $y^{[k]}$
в (\ref{eq9.4}) и (\ref{eq9.5}) совпадают построенными в~\S~\ref{s7}. Тогда в силу следствия~2, условия (\ref{eq9.4}) совпадают с (\ref{eq1.5}). Таким образом, для применения теоремы~\ref{th1} остается
показать, что условия (\ref{eq9.5}) могут быть представлены в виде (\ref{eq1.10}). Действительно, в силу (\ref{eq9.6}) условия (\ref{eq9.5}) можно считать нормированными следующим образом:
\begin{equation}\label{eq9.7}
U_j(y):=\sum_{k=0}^{N-1}\alpha_{j,k}y^{[k]}(0) +\sum_{k=0}^{\sigma_j}\beta_{j,k}y^{[k]}(1)=0, \quad \beta_{j,\sigma_j}=1, \quad j=\overline{l+1,N},
\end{equation}
причем выполняется (\ref{eq1.9}).  Подставляя (\ref{eq7.21}) в (\ref{eq9.7}), будем иметь
\begin{equation}\label{eq9.8}
U_j(y)=\sum_{k=0}^{\sigma_j}\Big(\beta_{j,k}y^{\langle k\rangle}(1) +\int_0^x y^{\langle k\rangle}(t)\,d\Phi_{j,k}(t)\Big) +\sum_{k=0}^{N-1}\alpha_{j,k}y^{\langle k\rangle}(0),
\end{equation}
где
$$
\Phi_{j,k}(t)=\sum_{\nu=k}^{\sigma_j}\beta_{j,\nu} \int_0^t R_{\nu,k}(1,\tau)\,d\tau, \quad k=\overline{0,\sigma_j}, \quad j=\overline{l+1,N}.
$$
С другой стороны, используя непосредственно определение (\ref{eq1.7}), а также формулу (\ref{eq6.7}), нетрудно получить представления
$$
y^{\langle k\rangle}(t)=\sum_{\nu=k}^{\sigma_j-1}y^{\langle \nu\rangle}(0)\frac{t^{\nu-k}}{(\nu-k)!}+ J^{\sigma_j-k}y^{\langle \sigma_j\rangle},
$$
если $\sigma_j<n$ и $k=\overline{0,\sigma_j-1}$ либо $\sigma_j>n$ и $k=\overline{n,\sigma_j-1},$ а также
$$
y^{\langle k\rangle}(t)=\sum_{\nu=k}^{n-1}y^{\langle \nu\rangle}(0)\frac{t^{\nu-k}}{(\nu-k)!} + J^{n-k}B^{-1}\Big(\sum_{\nu=n}^{\sigma_j-1}y^{\langle \nu\rangle}(0)\frac{t^{\nu-k}}{(\nu-k)!}
+J^{\sigma_j-n}y^{\langle \sigma_j\rangle} -Cy\Big),
$$
если $\sigma_j\ge n$ и $k=\overline{0,n-1}.$ Подставляя их в (\ref{eq9.8}) с учетом представления
$$
B^{-1}f(x)=f(x)+\int_0^x R(x,t)f(t)\,dt
$$
с непрерывным ядром $R(x,t),$ получаем (\ref{eq1.10}), где $\Phi_j\in AC[0,1].$ $\hfill\Box$
\end{proof}

\begin{remark}\label{rem8}
Доказательство теоремы~\ref{th10} фактически проведено для случая, когда вместо (\ref{eq9.5}) наложены более общие условия вида
$$
\sum_{k=0}^{N-1}\alpha_{j,k}y^{[k]}(0) +\sum_{k=0}^{\sigma_j}\Big(\beta_{j,k}y^{[k]}(1)+\int_0^1 y^{[k]}(t)\,d\varphi_{j,k}(t)\Big)=0, \quad \varphi_{j,k}\in BV[0,1].
$$
В самом деле, подставляя в них (\ref{eq7.21}), получаем (\ref{eq9.8}) с $\Phi_{j,k}\in BV[0,1].$
\end{remark}

\begin{remark}\label{rem9}
Выражение в (\ref{eq2.1}) не является частным случаем выражения $\ell_Ny,$ даже когда коэффициенты последнего подчиняются только общим условиям в разделе~\ref{s3.4}. Однако выражение в (\ref{eq2.1}) также можно
записать в виде (\ref{eq1.2})--(\ref{eq1.4}) (при $m=1$ и $n=N-1)$ с непрерывным ядром $K(x,t),$ равным нулю на диагонали, и представить краевые условия (\ref{eq2.2}) в виде (\ref{eq1.5}), (\ref{eq1.10}).
Поэтому утверждение теоремы~\ref{th2} тоже следует из теоремы~\ref{th1}.
\end{remark}

\end{fulltext}

\end{document}